\documentclass[11pt,A4]{article}

\usepackage{times}
\usepackage{amsmath}
\usepackage{amssymb}
\usepackage{amsfonts}
\usepackage{mathrsfs}
\usepackage{epsfig}
\usepackage{amsfonts}
\usepackage{amsthm}
\usepackage{graphics}
\usepackage{graphicx}
\usepackage{color}
\usepackage{eucal}

\definecolor{red}{rgb}{0.85,0,0}
\definecolor{green}{rgb}{0,0.6,0}
\definecolor{blue}{rgb}{0,0,0.85}

\renewcommand{\baselinestretch}{1.5}

\def\l{\mathcal{l}}
\def\T{\mathcal{T}}
\def\U{\mathcal{U}}

\def\abs#1{\left|#1\right|}

\def\~{\widetilde}

\def\0{\boldsymbol{0}}

\def\a{\bold{a}}
\def\e{\bold{e}}
\def\f{\bold{f}}
\def\w{\bold{w}}
\def\d{\bold{d}}
\def\l{\bold{l}}
\def\m{\bold{m}}
\def\n{\bold{n}}
\def\h{\bold{h}}
\def\u{\bold{u}}
\def\v{\bold{v}}
\def\r{\bold{r}}
\def\I{\bold{I}}
\def\J{\bold{J}}
\def\K{\bold{K}}

\def\vM{\boldsymbol{M}}
\def\vK{\boldsymbol{K}}

\def\vR{\boldsymbol{R}}
\def\vI{\boldsymbol{I}}
\def\vS{\boldsymbol{S}}
\def\vu{\boldsymbol{u}}
\def\vv{\boldsymbol{v}}
\def\vf{\boldsymbol{f}}
\def\vl{\boldsymbol{l}}

\def\vq{\boldsymbol{q}}

\def\vg{\boldsymbol{g}}
\def\pt{\partial_t}
\def\ps{\partial_s}
\def\half{\frac{1}{2}}

\title{Three Dimensional Nonlinear Dynamics of Slender Structures:\\
Cosserat Rod Element Approach}

\author{D. Q. Cao\footnote{
Corresponding Author. Email: d.cao@lancaster.ac.uk (DQ Cao)} ,
Dongsheng Liu and Charles H.-T. Wang\footnote{ Permanent address
starting 1 Jan 2005: School of Engineering and Physical Sciences,
University of Aberdeen, Aberdeen AB24 3UE, Scotland} }

\date{Department of Physics,
Lancaster University, Lancaster LA1 4YB, UK}

\begin{document}

\maketitle \vspace{-12pt}

\noindent{\bf Abstract}

In this paper, the modelling strategy of a Cosserat rod element
(CRE) is addressed systematically for 3-dimensional dynamical
analysis of slender structures. We employ the exact nonlinear
kinematic relationships in the sense of Cosserat theory, and adopt
the Bernoulli hypothesis. For the sake of simplicity, the Kirchoff
constitutive relations are adopted to provide an adequate
description of elastic properties in terms of a few elastic
moduli. A deformed configuration of the rod is described by the
displacement vector of the deformed centroid curves and an
orthonormal moving frame, rigidly attached to the cross-section of
the rod. The position of the moving frame relative to the inertial
frame is specified by the rotation matrix, parametrized by a
rotational vector. The approximate solutions of the nonlinear
partial differential equations of motion in quasi-static sense are
chosen as the shape functions with up to third order nonlinear
terms of generic nodal displacements. Based on the Lagrangian
constructed by the Cosserat kinetic energy and strain energy
expressions, the principle of virtual work is employed to derive
the ordinary differential equations of motion with third order
nonlinear generic nodal displacements. A simple example is
presented to illustrate the use of the formulation developed here
to obtain the lower order nonlinear ordinary differential
equations of motion of a given structure. The corresponding
nonlinear dynamical responses of the structures have been
presented through numerical simulations by Matlab software.


\newpage

\noindent \section{Introduction}

Three-dimensional slender structures undergoing large
displacements and rotations are often encountered in various
engineering systems such as vehicles, space structures, robotics,
aircrafts, and microelectronic mechanical systems. Clearly, these
systems consist of a set of interconnected components which may be
rigid or deformable. For example, a typical MEMS device may
consist of relatively heavy load bodies and thin springlike
supports. For such a system, each heavy body can be assumed to be
a rigid body and each springlike component can be described as a
deformable body. Since each of interconnected components of such a
system may undergo large displacements and/or rotations, an
effective modelling strategy that addresses very well to the
strongly nonlinear dynamic behavior is crucial in estimating
system performance and guiding the reliability verification
process.

Nonlinear finite element method provides a general approach to
structural modelling of multibody systems that consist of
interconnected {\it rigid} and {\it deformable} components. A
number of papers has recently been published, presenting new
concepts and new algorithms for modelling highly flexible spatial
frame structures \cite{Argyris, Cardona, Dutta}. An overview and
comprehensive treatment of this topic can be found, for instance,
in \cite{Shabana, Belytschko}. The Cosserat approach, that can
accommodate a good approximation the nonlinear behavior of complex
structures composed of materials with different constitutive
properties, variable geometry and damping characteristics
\cite{Green, Antman, Tucker, AntmanMarlow}, has been utilized to
develop finite element formulations for deformable bodies. The
finite element approach based on the Cosserat theory
(geometrically exact finite-strain beam theory) is usually
attributed to Reissner \cite{Reissner} and Simo \cite{Simo}. Simo
\cite{Simo} has discussed a convenient parameterization of the rod
model developed by Antman \cite{Antman72}, and Simo and Vu-Quoc
\cite{SimoVu} have considered the associated finite element
formulation. The computational procedure in \cite{SimoVu} uses a
variational formulation of the equations of motion and an
expansion of the kinematic quantities in terms of shape functions
and nodal values. Many modern finite element developers of the
three-dimensional beam theories, e.g. Jelenic and Saje
\cite{Jelenic}, Smolenski \cite{Smolenski}, and Zupan and Saje
\cite{Zupan} based their approach on the geometrically exact beam
theory. Another approach based on a system of Cosserat-type bodies
can be traced back to the work of Wozniak \cite{Wozniak}.
Homogeneously deformable bodies have been analyzed as pseudo-rigid
bodies \cite{Cohen} and Cosserat points \cite{Rubin1}. The theory
of a Cosserat point is a special continuum theory that models
deformation of a small structure that is essentially a point
surrounded by some small but finite region. The numerical
procedure based on the theory of a Cosserat point proposed in
\cite{Rubin1, Rubin2} has been used to study the dynamics of
spherically symmetric problems in \cite{Rubin3}. Recently, the
theory of a Cosserat point has been generalized to model a fully
nonlinear finite element for the numerical solution of 3-D dynamic
problems of elastic beams \cite{Rubin4}.

However, in practice the use of FEM codes to simulate complex
multibody systems such as MEMS devices is prohibitively
cumbersome, expensive, and time consuming. On the other hand, FEM
models use numerous variables to describe the device state. This
may lead the process of mapping the design space complex and the
relationship between each of these variables and the overall
device performance is not clear to designers. Recently, component
level modelling methods, which contain a library of parameterised
behavioural models for frequently used MEMS components
\cite{Lorenz, Mukherjee}, have been developed. In \cite{Lorenz,
Mukherjee}, every component is described as a single element in
contrast to FE models where the component is normally discretised
into many elements. Consequently, lower degree models are
established and the simulation time can be greatly reduced. The
mechanical behaviour of the components, however, is often modelled
using basic models, containing e.g. linear stiffness relationships
and/or approximations of basic nonlinearities.

Recently, motivated by the developments in MEMS modelling, the
Cosserat theory has been employed to develop a novel modelling
strategy that addresses very well to the practical needs for rapid
modelling of slender structures such as the springlike components
in MEMS, see Wang {\it et al.} \cite{WLR}. This modelling strategy
has been successfully used to investigate the non-ideal properties
of typical MEMS beams \cite{Gould}. In the sense of Cosserat
theory, the motion of rods in three-dimensional space can be
demonstrated by behaviors of a reference curve and three
perpendicular unit vectors (directors). Consequently, the
equations of motion are nonlinear partial differential equations,
which are functions of time and one space variable. For static
problems, however, the equations become nonlinear ordinary
differential equations, which can be solved approximately using
standard techniques like the perturbation method to satisfy
boundary conditions. In contrast, for dynamical problems, it is
necessary to introduce a numerical procedure which discretizes the
equations. In the strategy for modelling of a Cosserat rod element
\cite{WLR}, the basic kinematic quantities are the position of a
point on the Cosserat curve and an orthogonal transformation that
define the rotation of an orthogonal triad attached to the
cross-section at each point of the Cosserat curve. This enables
description of a rod using nonlinear ordinary differential
equations in terms of the generic nodal displacements of a CRE.

As an initial consideration, the modelling strategy in \cite{WLR}
is developed for 2-D case. In this paper, the modelling strategy
of CRE is addressed systematically for the 3-D problems. The
fundamental problem of any finite element formulation is the
choice of the shape functions. The approximate solutions of the
nonlinear equations of motion in quasi-static sense are chosen as
the shape functions with up to third order nonlinear terms of
generic nodal displacements. In three dimensions, the nonlinear
differential equations cannot be integrated in a closed form even
in the static sense, therefore the perturbation method is employed
here to solve the system approximately. For the sake of
simplicity, the Kirchoff constitutive relations are adopted to
provide an adequate description of elastic properties in terms of
a few elastic moduli. Based on the Lagrangian constructed by the
Cosserat kinetic energy and strain energy expressions, the
principle of virtual work is used to derive the ordinary
differential equations of motion with third order nonlinear
generic nodal displacements. The essential features and novel
aspects of the present formulation for CREs are briefly summarized
below.

\begin{enumerate}
\item The shape functions for CREs are derived from the
differential equations governing the flexural-flexural-torsional
motion of extensional rods, taking into account all the geometric
nonlinearities in the system. Consequently, the higher accuracy of
the dynamic responses can be achieved by dividing the rod into a
few elements which is much less than the traditional finite
element methods in which the interpolation functions are usually
extremely simple functions such as low order polynomials.

\item The mathematical simplicity when formulating deformable
bodies enables more convenient for modelling the multibody systems
that consist of interconnected {\it rigid} and {\it deformable}
components.

\item The resulting nonlinear ordinary differential equations with
lower degree-of-freedom are typically easy to simulate or
integrate into system-level simulations.
\end{enumerate}

An outline of the main contents of this paper is as follows. We
begin in section 2 by introducing the basic definitions and
kinematic assumptions on the nonlinear elastic rods that can
suffer flexure, extension, torsion, and shear. The rotational
vector that is free both of singularities and constraints is
employed as a parametrization to specify the deformed
configuration space. We limited our attention here to the
modelling of Cosserat rod elements in which the small effect of
shear will be neglected. The governing equations of motion and the
Kirchoff constitutive relations are presented in Section 3. The
straightforward perturbation method is employed in Section 4 to
solve the corresponding static problem. The approximate solutions
obtained are subsequently used as shape functions of Cosserat rod
elements. In section 5, Lagrangian approach is employed to
formulate the nonlinear ordinary differential equations of motion
of Cosserat rod elements. In terms of the shape functions derived
in Section 4, the Lagrangian is constructed by the Cosserat
kinetic energy and strain energy expressions, and the virtual work
done by external point loads and distributed loads is discussed. A
simple example is presented in section 7 to illustrate the use of
the formulation developed here to obtain the lower order nonlinear
ordinary differential equations of motion of a given structure.
The corresponding nonlinear dynamical responses of the structure
have been presented through numerical simulations by Matlab
software.

The following conventions and nomenclature will be used through
out this paper. Vectors, which are elements of Euclidean 3-space
$\mathscr{R}^3$, are denoted by lowercase, bold-face symbols,
e.g., $\u$, $\v$; vector-valued functions are denoted by
lowercase, italic, bold-face symbols, e.g., $\vu$, $\vv$; tensors
are denoted by upper-case, bold-face symbols, e.g., $\I$, $\J$;
matrices are denoted by upper-case, italic, bold-face symbols,
e.g., $\vM$, $\vK$. The three vectors $\{\e_1, \e_2, \e_3\}$ are
assumed to form a fixed right-handed orthogonal basis. The
summation convention for repeated indices is used. The symbols
$\pt$ and $\ps$ denote differentiation with respect to time $t$
and arc-length parameter $s$, respectively. The symbols $(\,\dot{}
\,)$ and $(\, ' \,)$ denote differentiation with respect to
dimensionless time parameter $\,\tau$ and dimensionless length
parameter $\,\sigma$, respectively.

\section{Kinematical Preliminaries}

\subsection{Basic definitions and kinematic assumptions}

\noindent Adopt Cartesian coordinates $(x,y,z)$ in inertial basis
$(\e_1,\e_2,\e_3)$ with Newtonian time $t$. According to the
Bernoulli hypothesis the plane cross-sections suffer only rigid
rotation during deformation and remain plane after deformation and
preserve their shape and area. For the sake of convenience, we
introduce the following definitions: (1) the reference
configuration, where the geometrical and mechanical variables of
the rod, including the loading, are known; (2) an arbitrary
deformed configuration, where only the loading is prescribed,
while the remaining variables are unknown.

\begin{figure}[ht]
\centerline{
\includegraphics[height=6.5cm]{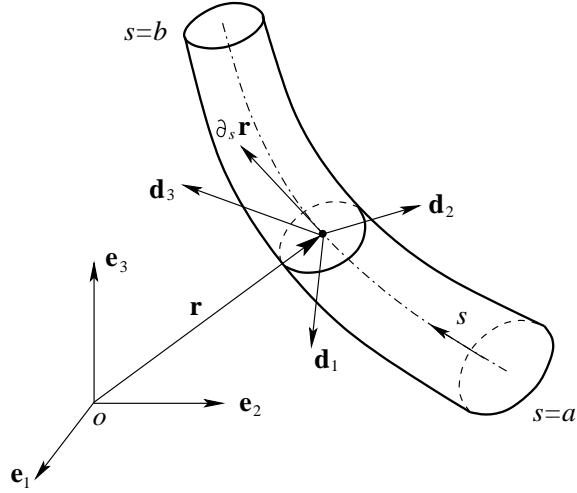}}
\caption{A simple Cosserat model} \label{figCmodel}
\end{figure}

It is therefore convenient to introduce an orthonormal basis
$\d_i(s,t), (i=1,2,3)$ of a cross-section at $s$, termed the
moving basis, such that $\d_3$ is normal to the rotated
cross-section, and $\d_1$ and $\d_2$ lie in the plane of the
rotated cross-section. The motion of a rod segment can be modelled
as a Cosserat rod whose configuration is described by its neutral
axis $\r(s,t)$ (Cosserat curve) and 3 orthogonal unit vectors
$\d_i(s,t), (i=1,2,3)$ (Cosserat directors) as shown in Figure
\ref{figCmodel}.

At any time, $\r$ describes the axis of the rod whose
cross-section orientations are determined by $\d_i$ such that $\ps
\r \cdot \d_3>0$. This condition implies that (i) the local ration
of deformed length to reference length of the axis cannot be
reduced to zero since $|\ps \r|>0$, and (ii) a typical
cross-section ($s=s_0$) cannot undergo a total shear in which the
plane determined by $\d_1$ and $\d_2$ is tangent to the curve
$\r(\cdot,t)$ at $\r(s_0,t)$ [10]. In general, as a result of
shear deformations of the rod, the cross-sections are not
perpendicular to the line of centroids.

In an inertial Cartesian basis $\{ \e_1, \e_2, \e_3 \}$ we may
write
\begin{equation}
\r(s,t) = r_i(s,t)\,\e_i = x(s,t) \e_1 + y(s,t) \e_2 + z(s,t)
\e_3. \label{position}
\end{equation}

The motion involves both the velocity of the curve, $\pt \r(s,t)$,
and angular velocity of the cross-sections $\w(s,t)$ so that
\begin{align}
\pt \d_i(s,t) = \w(s,t) \times \d_i(s,t). \label{def01}
\end{align}
In a similar manner the strains of the Cosserat rod are classified
into ``linear strain'' vector $\v(s,t) = \ps \r(s,t)$ and
``angular strain'' vector $\u(s,t)$ so that
\begin{align}
\ps \d_i(s,t) = \u(s,t) \times \d_i(s,t). \label{def02}
\end{align}

It follows from the definition (\ref{def01}) that
\begin{align}
\d_i \times \pt \d_i =  \d_i \times (\w \times \d_i)  = \w (\d_i
\cdot \d_i) - \d_i(\d_i \cdot \w ) = 2 \w. \nonumber
\end{align}
Therefore,
\begin{align}
\w = \frac{1}{2}\; \d_i \times \pt \d_i. \label{relation1}
\end{align}
Similarly, from the definition (\ref{def02}) we have
\begin{align}
\u = \frac{1}{2}\; \d_i \times \ps \d_i. \label{relation2}
\end{align}
Since the basis $\{\d_1, \d_2, \d_3\}$ is natural for the
intrinsic description of deformation, we decompose relevant vector
valued functions with respect to it:
\begin{align}
\v(s,t) =  v_i(s,t) \d_i(s,t), \qquad \u(s,t) = u_i \d_i(s,t),
\qquad \w(s,t) =  w_i \d_i(s,t). \label{vuw}
\end{align}

\subsection{Parametrization of the rotation matrix}

There is a number of choices for the parametrization of rotation
matrix, for example, the Euler angles, the quaternion parameters,
and the rotational vector being the most usual \cite{Stuelpnagel}.
Here, we employ the rotational vector that is free both of
singularities and constraints. Because of the orthogonality the
rotation matrix is a proper orthogonal matrix in SO(3), its nine
components can be expressed by only three independent parameters.
Denote $\vS$ the spin matrix of a vector $\a = a_i \e_i$ as
\begin{align}
\vS(\a) = \left[ \begin{array}{ccc} 0 & \,-a_3\, &
a_2 \\ a_3 & 0 & \, -a_1\, \\ \, -a_2 \, & a_1 & 0 \end{array}
\right]. \label{def03}
\end{align}
Then, the rotation matrix $\vR$ is determined by the expression
\cite{Stuelpnagel}
\begin{align}
\vR(\boldsymbol{\phi}) = \vI + \frac{\sin \phi}{\phi}
\vS(\boldsymbol{\phi}) + \frac{1 - \cos \phi}{\phi^2}
\vS^2(\boldsymbol{\phi}), \label{def04}
\end{align}
where $\boldsymbol{\phi} = \phi_i \e_i$ is the rotational vector,
$\vS(\boldsymbol{\phi})$ is the spin matrix of $\boldsymbol{\phi}$
defined by (\ref{def03}), and $\phi = (\phi_1^2 + \phi_2^2 +
\phi_3^2)^{1/2}$ is the rotational norm or the length of the
rotational vector. An expansion of trigonometric functions in Eq.
(\ref{def04}) in MacLaurin's series yields
\begin{align}
\vR = \vI + \vS + \frac{1}{2!} \vS^2 + \frac{1}{3!} \vS^3 + \cdots
+ \frac{1}{n!} \vS^n + \cdots = \exp \vS. \label{Mcseries}
\end{align}
Thus, the rotation matrix may alternatively be expressed by an
exponential map, the exponentiation of the spin matrix associated
with the rotational vector. Note that, as a consequence of the
exponentiation of the spin matrix $\vS(\boldsymbol{\phi})$ being
equal to $\vR(\boldsymbol{\phi})\in$ SO(3), the spin matrix
$\vS(\boldsymbol{\phi})$ belongs to Lie algebra so(3) associated
with the Lie group SO(3) \cite{Jones}.

Conversely, taking a given orthogonal matrix $\vR$ as a rotation
matrix, the associated rotation vector $\boldsymbol{\phi}$ can be
derived from (\ref{def03}) and (\ref{def04}). The rotational norm
$\phi$ can be calculated by
\begin{align}
\phi=\cos ^{-1} \frac{Tr(\vR) - 1}{2}.
\end{align}
By taking the matrix logarithm of $\vR$ we can obtain the
skew-symmetric matrix $\vS$ as following.
\begin{align}
\vS = \log \vR = \frac{\phi}{2 \sin \phi}(\vR-\vR^T).
\end{align} \label{loga}
Therefore $\boldsymbol{ \phi} = \phi_i \e_i$ with $\phi_1= -
S_{23},\, \phi_2 = S_{13}$, and $\phi_3 = - S_{12}$.

In terms of the rotational vector $\boldsymbol{\phi}$, Eqs.
(\ref{def03}) and (\ref{def04}) give the exact value of the
current rotation matrix. Using truncated MacLaurin's series of
various order in Eq. (\ref{Mcseries}), approximate values of the
rotation matrix are obtained and corresponding simplified theories
can be derived. For example, a so called first order theory is
obtained if small rotations are assumed so that the quadratic and
higher order terms in Eq. (\ref{Mcseries}) may be neglected.

\subsection{Specifications for the deformed configuration space}

The position vector $\r(s,t)$ defined by (\ref{position}) is an
element of Euclidean vector space $\mathscr{R}^3$. The orientation
of the moving basis is represented by the rotation matrix, which
is en element of the Lie group SO(3). Accordingly, the set of all
possible configurations of the rod is defined by
\begin{align}
\mathscr{C} = \{ (\r, \vR) \vert\, \r: s \to \mathscr{R}^3, \,
\vR: s \to \textnormal{SO}(3) \}.
\end{align}
This set is referred to as the deformed configuration space. The
quantities $\r$ and $\vR$ are termed the kinematic quantities of
the rod. Since the rotation matrix is related to the three
parameters, the components of the rotational vector
$\boldsymbol{\phi}$, the Lie group SO(3) of rotation matrices is
three-parametric, i.e. it may be viewed as being a 3-D nonlinear
differentiable manifold.

For a typical slender rod such as the components in MEMS, the
effect of shearing deformation can be negligible, the
cross-section of the rod is therefore assumed to be perpendicular
to the tangent to the Cosserat curve, i.e.
\begin{align}
\v(s,t) = \ps \r(s,t) = \abs{\ps \r(s,t)} \d_3(s,t).
\label{nosh01}
\end{align}
In this case, we write
\begin{align}
\d_3(s,t) = \frac{\ps\r(s,t)}{\abs{\ps \r(s,t)}} \triangleq
\nu_1(s,t) \e_1 + \nu_2(s,t) \e_2 + \nu_3(s,t) \e_3 \label{nosh02}
\end{align}
with
\begin{align}
\nu_1^2(s,t) + \nu_2^2(s,t) + \nu_3^2(s,t) = 1.
\end{align}
where $\nu_1, \nu_2$ and $\nu_3$ can be written as
\begin{align}
\nu_1(s,t) = \frac{\ps x(s,t)}{\abs{\ps\r(s,t)}}, \quad \nu_2(s,t)
= \frac{\ps y(s,t)}{\abs{\ps\r(s,t)}}, \quad \textnormal{and}
\quad \nu_3(s,t) = \frac{\ps z(s,t)}{\abs{\ps\r(s,t)}},
\label{mu123}
\end{align}
by differentiating the position vector $\r(s,t)$ defined in
(\ref{position}) with respect to $s$.

We assume that the directors $\{\d_1,\d_2,\d_3\}$  can be obtained
by the following way. First of all, we rotate directors
$\{\e_1,\e_2,\e_3\}$ about $\e_3$ with an angle $\varphi$ to
obtain the directors $\{\tilde \d_1,\, \tilde \d_2,\, \e_3\}$.
Then, rotation matrix $\vR_a$ associated with the rotational
vector $\boldsymbol{ \phi_a} = \varphi \, \e_3$ can be written as
\begin{align}
\vR_a = \vR(\boldsymbol{ \phi_a}) = \left[
\begin{array}{ccc}
\; \cos \varphi \; & \; - \sin \varphi\; & \;0 \; \\
\sin \varphi & \cos \varphi & 0 \\ 0 & 0 & 1  \end{array} \right].
\label{rotate1}
\end{align}

Next, we introduce a rotational vector
\begin{align}
\boldsymbol{ \phi_b} = - \frac{\sin^{-1}\sqrt{\nu_1^2 + \nu_2^2}
}{\sqrt{\nu_1^2+\nu_2^2}} \, \nu_2 \, \tilde \d_1 +
\frac{\sin^{-1} \sqrt{\nu_1^2 + \nu_2^2}}{\sqrt{\nu_1^2+\nu_2^2}}
\, \nu_1 \, \tilde \d_2 \nonumber
\end{align}
which rotates the vectors $\{\tilde \d_1,\, \tilde \d_2, \,\e_3\}$
to $\{\d_1,\, \d_2, \,\d_3\}$. Here, we assume that $\nu_1^2 +
\nu_2^2 \ne 0$. Other wise $\d_3 = \e_3$, this rotating procedure
can be omitted. Let $\vR_b$ be the corresponding rotation matrix
associated with the rotational vector $\boldsymbol{ \phi_b}$. Then
\begin{align}
\vR_b = \vR(\boldsymbol{ \phi_b}) = \left[ \begin{array}{ccc}
\frac{\nu_1^2\nu_3 + \nu_2^2}{\nu_1^2+\nu_2^2} & \;
\frac{\nu_1 \nu_2(\nu_3-1)}{\nu_1^2+\nu_2^2}\; & \; \nu_1 \\
\; \frac{\nu_1\nu_2(\nu_3-1)}{\nu_1^2+\nu_2^2} \; &
 \frac{\nu_2^2\nu_3+\nu_1^2}{\nu_1^2+\nu_2^2} & \nu_2\\
-\nu_1 & \nu_2 & \nu_3
\end{array}\right]\; .
\end{align}
Consequently, the moving directors are obtained as:
\begin{align} \label{md1}
\d_1 = & \left( \frac{(\nu_1^2\nu_3+\nu_2^2)\cos
\varphi}{\nu_1^2+\nu_2^2} + \frac{\nu_1\nu_2(\nu_3-1)\sin
\varphi}{\nu_1^2 + \nu_2^2} \right) \, \e_1 \nonumber \\[6pt]
& + \left( \frac{(\nu_2^2\nu_3+\nu_1^2)\sin \varphi
}{\nu_1^2+\nu_2^2} + \frac{\nu_1\nu_2(\nu_3-1)\cos
\varphi}{\nu_1^2+\nu_2^2} \right) \, \e_2
- (\nu_1\cos \varphi + \nu_2\sin \varphi)\, \e_3, \\[6pt]
\d_2 = & \left( - \frac{(\nu_1^2\nu_3+\nu_2^2) \sin \varphi}
{\nu_1^2+\nu_2^2} + \frac{\nu_1\nu_2(\nu_3-1) \cos \varphi}
{\nu_1^2 + \nu_2^2} \right) \, \e_1 \nonumber\\[6pt]
& + \left( \frac{(\nu_2^2\nu_3 + \nu_1^2) \cos
\varphi}{\nu_1^2+\nu_2^2} - \frac{\nu_1\nu_2(\nu_3-1) \sin
\varphi}{\nu_1^2+\nu_2^2} \right) \, \e_2
+(\nu_1 \sin \varphi - \nu_2 \cos \varphi) \, \e_3, \\[6pt]
\d_3 = & \, \nu_1\e_1 + \nu_2\e_2 + \nu_3 \e_3. \label{md3}
\end{align}
Obviously  $\varphi(s,t)$ is a variable related to torsion of the
rod. Expanding the directors in polynomials about
$\nu_1,\nu_2,\phi$ and reserving the terms up to third order, we
have
\begin{align} \label{di1}
\d_1(s,t) \approx & \left( 1 - \frac{1}{2} \varphi^2(s,t) -
\frac{1}{2} \nu_1^2(s,t) - \frac{1}{2} \nu_1(s,t) \nu_2(s,t)
\varphi(s,t)\right)\, \e_1 \nonumber\\[4pt]
& + \left( \varphi(s,t) - \frac{1}{2} \nu_1(s,t)\nu_2(s,t) -
\frac{1}{2} \nu_2^2(s,t) \varphi(s,t) - \frac{1}{6} \varphi^3(s,t)
\right) \, \e_2 \nonumber\\[4pt]
& + \left( -\nu_1(s,t) - \nu_2(s,t)\varphi(s,t) + \frac{1}{2}
\nu_1(s,t) \varphi^2(s,t) \right)\, \e_3 \\[4pt]
\d_2(s,t) \approx & \left( - \varphi(s,t) - \frac{1}{2}
\nu_1(s,t)\nu_2(s,t) + \frac{1}{2} \nu_1^2(s,t)\varphi(s,t) +
\frac{1}{6}\phi^3(s,t) \right) \,\e_1 \nonumber\\[4pt]
& + \left( 1 - \frac{1}{2} \varphi^2(s,t) -
\frac{1}{2}\nu_2^2(s,t) + \frac{1}{2}\nu_1(s,t)\nu_2(s,t)
\varphi(s,t) \right)\, \e_2 \nonumber\\[4pt] & + \left(
-\nu_2(s,t) + \nu_1(s,t)\phi(s,t) + \frac{1}{2}
\nu_2(s,t) \varphi^2(s,t) \right) \e_3 \\[4pt]
\label{di3} \d_3(s,t) \approx & \, \nu_1(s,t) \e_1 +
\nu_2(s,t)\e_2 + \left( 1 - \frac{1}{2}\nu_1^2(s,t) -
\frac{1}{2}\nu_2^2(s,t) \right) \, \e_3
\end{align}

For the sake of convenience to describe the displacements and
rotations in the inertia frame, we regard directors $\d_i(s,t)$
$(i=1,2,3)$ as those obtained by rotating inertial frame
$\{\e_1,\e_2,\e_3\}$ with a rotation vector
\begin{align}
\boldsymbol{\phi} = \phi_{x}(s,t)\e_1 + \phi_{y}(s,t)\e_2 +
\phi_{z}(s,t)\e_3. \label{rphi}
\end{align}
Now, based on the relations (\ref{md1})--(\ref{md3}), utilizing
the inverse procedure mentioned in Section 2.2, the rotational
norm $\phi$ and the spin matrix associated with the rotation
vector $\boldsymbol{\phi}$ in (\ref{rphi}) can be derived from
\begin{align}
\phi = \cos ^{-1} \frac{Tr(\vR_b \vR_a) - 1}{2}.
\end{align}
and
\begin{align}
\vS = \log (\vR_b \vR_a) = \frac{\phi}{2 \sin \phi}(\vR_b \vR_a -
\vR_a^T \vR_b^T). \label{loga1}
\end{align}

Consequently, the approximate relations between $(\phi_x, \phi_y,
\phi_z)$ and $(\nu_1, \nu_2, \varphi)$, up to third order, are
obtained as
\begin{align} \left\{ \begin{array}{l} \displaystyle
\phi_{x}(s,t) = - \nu_2(s,t) + \frac{1}{2}\,
\varphi(s,t)\nu_1(s,t) -\frac{1}{6}\, \left( \nu_1^2(s,t) +
\nu^2_2(s,t)
- \frac{1}{2}\, \varphi^2(s,t) \right) \nu_2(s,t),  \\[10pt]
\displaystyle \phi_{y}(s,t) = \nu_1(s,t) + \frac{1}{2}\,
\varphi(s,t)\nu_2(s,t) + \frac{1}{6}\, \left( \nu_1^2(s,t) +
\nu^2_2(s,t)
- \frac{1}{2}\, \varphi^2(s,t) \right) \nu_1(s,t), \\[10pt]
\displaystyle \phi_{z}(s,t) = \varphi(s,t) - \frac{1}{12}\, \left(
\nu_1^2(s,t) + \nu_2^2(s,t) \right) \varphi(s,t),
\end{array} \right. \label{rota}
\end{align}
or equivalently,
\begin{align} \left\{ \begin{array}{l} \displaystyle
\nu_1(s,t) = \phi_{y}(s,t) + \frac{1}{2}\,
\phi_{x}(s,t)\phi_{z}(s,t) - \frac{1}{6}\, \left( \phi^2_x(s,t) +
\phi_y^2(s,t) + \phi^3_z(s,t)
\right) \phi_y(s,t), \\[10pt] \displaystyle
\nu_2(s,t) = -\phi_{x}(s,t) + \frac{1}{2}\, \phi_{y}(s,t)
\phi_{z}(s,t) + \frac{1}{6}\, \left( \phi^2_x(s,t) + \phi^2_y(s,t)
+ \phi_z^2(s,t) \right) \phi_x(s,t), \\[10pt] \displaystyle
\varphi(s,t) = \phi_{z}(s,t) + \frac{1}{12}\, \left( \phi_x^2(s,t)
+ \phi_y^2(s,t) \right) \phi_z(s,t).
\end{array} \right. \label{ang}
\end{align}

These relations are very useful in solving the static problem and
will be used below to derive the shape functions for CRD.

\section{The Governing Equations of Motion}

The dynamical evolution of the rod with density, $\rho(s)$, and
cross-section area, $A(s)$ is governed by the Newton's dynamical
laws:
\begin{align}\left\{ \begin{array}{ll}
& \rho(s) A(s) \partial_{tt} \r = \ps \n(s,t) + \f(s,t), \\
& \pt \h(s,t) = \ps \m(s,t) + \v(s,t) \times \n(s,t) + \l(s,t),
\end{array} \right. \label{dequation1}
\end{align}
where
\begin{align}
\n(s,t) = n_i(s,t) \d_i(s,t), \quad \m(s,t) = m_i(s,t) \d_i(s,t)
\end{align}
are the contact force and contact torque densities, respectively;
while
\begin{align}
\h(s,t) = h_i(s,t) \d_i(s,t)
\end{align}
denotes the angular momentum densities; $\f(s,t)$ and $\l(s,t)$
denote the prescribed external force and torque densities,
respectively.

The simplest constitutive model is based on the Kirchoff
constitutive relations which provide an adequate description of
elastic properties in terms of a few elastic moduli. One may
exploit the full versatility of the Cosserat model by generating
the Kirchoff constitutive relations to include viscoelasticity and
other damping, curved reference states with memory and effects to
prohibit total compression.

The contact forces, contact torques and the angular momentum are
given as
\begin{align}
\n = \K (\v - \d_3), \qquad \m = \J (\u), \qquad \h = \I (\w)
\label{ftm}
\end{align}
where, according to the Kirchoff constitutive relations, the
tensors $\K$, $\J$ and $\I$ are described as
\begin{align} \left\{ \begin{array}{ll}
 & \K(s,t) = \displaystyle K_{ii}(s,t)(\d_i(s,t) \otimes \d_i(s,t)),\\
 & \J(s,t) = \displaystyle \sum_{i,j=1}^2 J_{ij}(s,t)(\d_i(s,t) \otimes
 \d_j(s,t)) + J_{33}(s,t)(\d_3(s,t) \otimes \d_3(s,t)),\\
 & \I(s,t) = \displaystyle \sum_{i,j=1}^2 I_{ij}(s,t)(\d_i(s,t) \otimes
 \d_j(s,t)) + I_{33}(s,t)(\d_3(s,t) \otimes \d_3(s,t)).
\end{array} \right. \label{Kirchoff}
\end{align}
The corresponding components are given as
\begin{align} \left\{ \begin{array}{lll}
& K_{11} = K_{22} = G A(s), & \quad K_{33} = EA(s),\\
& J_{11} = \int_{A(s)} E \eta^2 \,dA, & \quad J_{22} =
\int_{A(s)} E \xi^2\, dA,\\
& J_{33} = \int_{A(s)} G (\xi^2+\eta^2) \,dA, & \quad J_{12} = -
J_{21} = \int_{A(s)} E \xi \eta\, dA, \\
& I_{11} = \int_{A(s)} \rho(s) \eta^2 \,dA, & \quad I_{22} =
\int_{A(s)} \rho(s) \xi^2\, dA,\\
& I_{33} = \int_{A(s)} \rho(s) (\xi^2+\eta^2) \,dA, & \quad I_{12}
= - I_{21} = \int_{A(s)} \rho(s) \xi \eta\, dA, \end{array}
\right. \label{tensor1}
\end{align}
where $E$ and $G$ are the Young's modulus of elasticity and shear
modulus respectively.

\section{Shape Functions for Cosserat Rod Elements}

For convenience, consider a uniform and initially straight rod
element of constant length $L$, supported in an arbitrary manner
at $s=a=0$ and $s=b=L$. It is assumed in the following that the
static equilibrium of the rod corresponds to the situation where
the directions of $\d_3$ and $\e_3$ are coincident with each other
and $\d_1$, $\d_2$ are parallel to $\e_1$, $\e_2$, respectively.
The principal axes are chosen to parallel $\e_1$, $\e_2$ and
$\e_3$. For the sake of simplicity, it will be assumed that the
axes along the directors $\d_1$, $\d_2$ and $\d_3$ are chosen to
be the principal axes of inertia of the cross section at $s$, and
centered at the cross section's center of mass. Then, for a
uniform rod with cross-section area $A(s)$, we have $J_{12} =
J_{21} = 0$,  $I_{12} = I_{21} = 0$ and
\begin{align} \left\{ \begin{array}{lll}
& K_{11} = K_{22} = G A(s), & \quad K_{33} = EA(s),\\
& J_{11} = E \int_{A(s)} \eta^2 \,dA, & \quad J_{22} = E \int_{A(s)} \xi^2\, dA,\\
& I_{11} = \rho \int_{A(s)} \eta^2 \,dA, & \quad I_{22} = \rho
\int_{A(s)} \xi^2\, dA,\\
& J_{33} = G \int_{A(s)} (\xi^2 + \eta^2)
\,dA = \frac{G}{E} (J_{11} + J_{22} ), & \\
& I_{33} = \rho \int_{A(s)} (\xi^2+\eta^2) \,dA = I_{11} + I_{22}.
& \end{array} \right. \label{tensor2}
\end{align}

Assume that the shape functions for a CRE satisfy the
corresponding static equations of (\ref{dequation1}), i.e.
\begin{align}
& \ps \n(s) = 0  \label{dequ} \\
& \ps \m(s) + \v(s) \times \n(s) = 0, \label{sequ}
\end{align}
where the contact force and contact torque densities are
\begin{align} \left\{ \begin{array}{l}
\, \n(s) = n_i(s) \d_i(s), \quad n_1 = K_{11} v_1, \quad n_2 =
K_{22} v_2, \quad n_3 = K_{33} ( v_3 - 1), \\
\, \m(s) = m_i(s) \d_i(s), \quad m_1 = J_{11} u_1, \quad m_2 =
J_{22} u_2, \quad m_3 = J_{33} u_3. \end{array} \right.
\end{align}
with $\u(s) = \frac{1}{2} \d_i(s) \times \ps \d_i(s)$, and
$\d_i(s)$ $(i=1,2,3)$ are given by (\ref{di1})--(\ref{di3}).

As mentioned in Section 2.3, for a typical slender rod as the
components in MEMS, the effect of shearing deformation can be
negligible, therefore  the cross-section of rod is assumed to be
perpendicular to the tangent to the Cosserat curve, i.e. the
strain vector $\v(s,t) = \abs{\ps \r(s)} \d_3(s)$ satisfies the
form (\ref{nosh01}). Thus, $v_1 = v_2 =0$ and $v_3 = \abs{\ps
\r(s)}$. Consequently, instead of $n_1 = K_{11} v_1$ and $n_2 =
K_{22} v_2$, the contact forces $n_1$ and $n_2$ follow from
(\ref{sequ})
\begin{align}
n_1 = \displaystyle \frac{- \ps m_2 - u_3 m_1 + u_1 m_3}{v_3},
\quad \textnormal{and} \quad n_2 = \displaystyle \frac{\ps m_1 -
u_3 m_2 + u_2 m_3}{v_3}. \label{nosh2}
\end{align}

As a prelude to expanding the nonlinear shape functions to a form
suitable for a perturbation analysis of the motion, it is useful
to introduce some natural scales to obtain a dimensionless
equation of motion. Introduce the dimensionless variables
\begin{align}
\sigma = \frac{s}{L_0}, \; \bar \r = \frac{\r}{L_0}, \; \bar x =
\frac{x}{L_0}, \; \bar y = \frac{y}{L_0}, \; \bar z =
\frac{z}{L_0}, \; \tau = \omega_0 t, \label{dless}
\end{align}
where $L_0$ and $\omega_0$ are the reference length and natural
frequency yet to be determined later, respectively.

Assume that the dimensionless generic nodal displacements
(boundary displacements and rotations) at $\sigma = 0$ and $\sigma
= L/L_0$ are
\begin{align}
\vq_a = \left[ \begin{array}{cccccc} \epsilon X_a & \epsilon Y_a &
\epsilon Z_a & \epsilon \Phi_{xa} & \epsilon \Phi_{ya} & \epsilon
\Phi_{za}
\end{array} \right]^T \label{staticnda}
\end{align}
and
\begin{align}
\vq_b = \left[ \begin{array}{cccccc} \epsilon X_b & \epsilon Y_b &
\epsilon Z_b & \epsilon \Phi_{xb} & \epsilon \Phi_{yb} & \epsilon
\Phi_{zb}
\end{array} \right]^T, \label{staticndb}
\end{align}
respectively. Substituting (\ref{staticnda}) and (\ref{staticndb})
into (\ref{position}), we obtain the boundary conditions for $\bar
x$, $\bar y$ and $\bar z$ as
\begin{align} \left\{ \begin{array}{lll}
\;\bar x(0) = \epsilon X_a, \quad  & \bar y(0) = \epsilon Y_a,
\quad & \bar z(0) = \epsilon Z_a,  \\
\;\bar x(l) = \epsilon X_b, \quad & \bar y(l) = \epsilon Y_b,
\quad & \bar z(l) = l + \epsilon Z_b, \end{array} \right.
\label{dipb01}
\end{align}
where $l = L/L_0$ is the dimensionless length of the rod element.
Substituting (\ref{staticnda}) and (\ref{staticndb}) into
(\ref{ang}), we obtain the boundary conditions for $\nu_1$,
$\nu_2$ and $\varphi$ as
\begin{align} \left\{ \begin{array}{l} \displaystyle
\; \nu_1(0) = \frac{\bar x'(0)}{\abs{\bar \r'(0)}} = \epsilon
\Phi_{ya} + \frac{1}{2} \epsilon^2 \Phi_{xa} \Phi_{za} -
\frac{1}{6} \epsilon^3 \left(
\Phi_{xa}^2 + \Phi_{ya}^2 + \Phi_{za}^2 \right) \Phi_{ya}, \\[9pt]
\displaystyle \; \nu_2(0) = \frac{\bar y'(0)}{\abs{\bar \r'(0)}} =
- \epsilon \Phi_{xa} + \frac{1}{2} \epsilon^2 \Phi_{ya} \Phi_{za}
+ \frac{1}{6} \epsilon^3 \left(
\Phi_{xa}^2 + \Phi_{ya}^2 + \Phi_{za}^2 \right) \Phi_{xa}, \\[9pt]
\displaystyle \; \varphi(0) = \epsilon \Phi_{za} + \frac{1}{12}
\epsilon^3 \left( \Phi_{xa}^2 + \Phi_{ya}^2 \right) \Phi_{za}
\\[9pt]
\displaystyle \; \nu_1(l) = \frac{\bar x'(l)}{\abs{\bar \r'(l)}} =
\epsilon \Phi_{yb} + \frac{1}{2} \epsilon^2 \Phi_{xb} \Phi_{zb} -
\frac{1}{6} \epsilon^3 \left(
\Phi_{xb}^2 + \Phi_{yb}^2 + \Phi_{zb}^2 \right) \Phi_{yb}, \\[9pt]
\displaystyle \; \nu_2(l) = \frac{\bar y'(l)}{\abs{\bar \r'(l)}} =
- \epsilon \Phi_{xb} + \frac{1}{2} \epsilon^2 \Phi_{yb} \Phi_{zb}
+ \frac{1}{6} \epsilon^3 \left(
\Phi_{xb}^2 + \Phi_{yb}^2 + \Phi_{zb}^2 \right) \Phi_{xb}, \\[9pt]
\displaystyle \; \varphi(l) = \epsilon \Phi_{zb} + \frac{1}{12}
\epsilon^3 \left( \Phi_{xb}^2 + \Phi_{yb}^2 \right) \Phi_{zb}
\end{array} \right. \label{rotb01}
\end{align}

Treating $\epsilon$ as a perturbation parameter which is the order
of the amplitude of the displacement and can be used as a crutch
in obtaining the approximate solution, the shape functions can be
obtained by solving the static equations (\ref{dequ}) and
(\ref{sequ}) with the corresponding boundary conditions
(\ref{dipb01}) and (\ref{rotb01}) and also the restrictions
(\ref{nosh2}) on the assumption of neglecting the effect of
shearing deformation. To do this, we seek a straightforward
expansion
\begin{align} \left\{ \begin{array}{l}
\, \bar x(\sigma) = \epsilon \hat x_1(\sigma) + \epsilon^2 \hat
x_2(\sigma)
+ \epsilon^3 \hat x_3(\sigma) + \cdots,  \\
\, \bar y(\sigma) = \epsilon \hat y_1(\sigma) + \epsilon^2 \hat
y_2(\sigma)
+ \epsilon^3 \hat y_3(\sigma) + \cdots,  \\
\, \bar z(\sigma) = \sigma + \epsilon \hat z_1(\sigma) +
\epsilon^2 \hat z_2(\sigma)
+ \epsilon^3 \hat z_3(\sigma) + \cdots,  \\
\, \varphi(\sigma) = \epsilon \hat \varphi_1(\sigma) + \epsilon^2
\hat \varphi_2(\sigma) + \epsilon^3 \hat \varphi_3(\sigma) +
\cdots . \end{array} \right. \label{shape}
\end{align}

Substituting (\ref{shape}) into (\ref{dequ}) and (\ref{sequ})
associated with (\ref{nosh2}) and, because $\bar x_i$, $\bar y_i$,
$\bar z_i$ and $\bar \varphi_i$ are independent of $\epsilon$, set
the coefficient of each power of $\epsilon$ equal to zero. This
leads to a set of linear ordinary differential equations which can
be solved using the Frobeniu's method \cite{Arfken} under the
corresponding boundary conditions (\ref{dipb01}) and
(\ref{rotb01}). The solving procedure has been implemented in a
MAPLE file \cite{LW}. Consequently, the approximate series
solutions are obtained and the 1st order ones are
\begin{align}  \left\{ \begin{array}{ll}
\hat x_1(\sigma) =  X_a & + \;\Phi_{ya} \sigma - (3X_a - 3X_b + 2l
\Phi_{ya} + l \Phi_{yb}) \displaystyle  \frac{\sigma^2}{l^2} \\
& + \;(2X_a - 2X_b + l \Phi_{ya} + l \Phi_{yb}) \displaystyle
\frac{\sigma^3}{l^3} \\
\hat y_1(\sigma) = Y_a & - \;\Phi_{xa} \sigma -(3Y_a - 3Y_b - 2l
\Phi_{xa} - l \Phi_{xb}) \displaystyle  \frac{\sigma^2}{l^2} \\
& + \;(2Y_a - 2Y_b - l \Phi_{xa} - l \Phi_{xb} )
\displaystyle  \frac{\sigma^3}{l^3} \\
\hat z_1(\sigma) =  Z_a & + \;(Z_b - Z_a) \displaystyle  \frac{\sigma}{l}\\
\hat \varphi_1(\sigma) = \Phi_{za} & + \;(\Phi_{zb} - \Phi_{za})
\displaystyle  \frac{\sigma}{l}. \end{array} \right.
\label{curve1}
\end{align}
To investigate deflections up to 3rd order nonlinearity in
$\epsilon$ it is adequate to adopt the truncated (\ref{shape}) to
$\epsilon^3$ order terms. The high order terms (up to third order)
which are polynomials of $\sigma$, can be easily solved using a
MAPLE programm \cite{LW}. For example, $\hat{x}_2(\sigma)= C_1
\sigma^5 + C_2 \sigma^4 + C_3 \sigma^3  + C_4 \sigma^2$ with
\begin{align} \label{curve3}
C_1 = \frac{K_{33}}{20l^4J_{22}} (Z_b - Z_a)(2X_a - 2X_b + l
\Phi_{ya} + l \Phi_{yb}).
\end{align}

Accordingly to the time-dependent, rod shape under the
quasi-static condition is specified with the (slowly) time-varying
nodal displacements and rotations.

\section{Equations of Motion for Cosserat Rod Elements}

\noindent In this section, the Lagrangian approach is employed to
formulate the ordinary differential equations of motion of
Cosserat rod elements. The generalized Hamilton's principle which,
in its most general form, is given by the variational statement
\begin{align}
\int_{t_1}^{t_2} \delta (\mathscr{T} - \mathscr{V}) dt +
\int_{t_1}^{t_2} \delta \mathscr{W} dt = 0, \label{Hamilton}
\end{align}
where $\mathscr{T}$ is the total kinetic energy of the system,
$\mathscr{V}$ is the potential energy of the system (including the
strain energy and the potential energy of conservative external
forces), $\delta(\,\cdot \,)$ represents the virtual displacement
(or variational) operator, and $\delta \mathscr{W}$ is the virtual
work done by nonconservative forces (including damping forces) and
external forces not accounted for in $\mathscr{V}$.

Assume that the time-varying dimensionless displacements at the
ends ($\sigma=a/L_0$ and $\sigma=b/L_0$) of the element model are
\begin{align}
\vq_a(\tau) = \left[ \begin{array}{cccccc} X_a(\tau) & Y_a(\tau) &
Z_a(\tau) & \Phi_{xa}(\tau) & \Phi_{ya}(\tau) & \Phi_{za}(\tau)
\end{array} \right]^T \label{dynamicnda}
\end{align}
and
\begin{align}
\vq_b(\tau) = \left[ \begin{array}{cccccc} X_b(\tau) & Y_b(\tau) &
Z_b(\tau) & \Phi_{xb}(\tau) & \Phi_{yb}(\tau) & \Phi_{zb}(\tau)
\end{array} \right]^T, \label{dynamicndb}
\end{align}
respectively. Then, the generalized displacement vector for the
element can be described by
\begin{align} \vq^e(\tau) = \left[
\begin{array}{cc} \vq_a^T(\tau) & \vq_b^T(\tau)
\end{array} \right]^T
\end{align}

Consistent with the kinematic and constitutive assumptions
described in Section 2 and Section 3 and the shape functions
derived in Section 4, the kinetic energy per unit length is
\begin{align} \label{k1}
\T & = \half \left\{ \rho A \,\pt\r \cdot \pt\r + \I(\w,\w)
\right\} = \half \left\{ \rho A \omega_0^2 L_0^2\,\dot {\bar \r}
\cdot \dot {\bar\r} + \omega_0^2 \I(\bar\w,\bar\w) \right\}
\end{align}
where $\rho$ and $A$ are the density of rod and the area of
cross-section of rod, respectively. According to (\ref{position})
and (\ref{relation1}), the velocity $\pt \r(s,t)$, and the angular
velocity of the cross-section can be derived as:
\begin{align}
\pt \r  = \pt x \e_1 + \pt y \e_2 + \pt z \e_3 = \omega_0 L_0
(\dot{\bar x} \e_1 + \dot{\bar y} \e_2 + \dot{\bar z} \e_3)  =
\omega_0 L_0 \dot{\bar{\r}}
\end{align}
and
\begin{align}
\w = \half \d_i \times \pt \d_i  = \half \omega_0 \d_i \times \dot
\d_i = \omega_0 \bar \w, \label{dimlessw}
\end{align}
respectively.

Under small strain conditions the strain energy per unit length of
rod can be expressed in terms of the strain vectors $\u$ and $\v$
as:
\begin{align}\label{p1}
\U & = \half \left\{ \J(\u,\u) + K_{33}(v_3 - 1)^2 \right \} =
\half \left\{\frac{1}{L_0^2} \J(\bar \u,\bar\u) + K_{33}(\bar v_3
- 1)^2 \right\}
\end{align}
where the strain vector is
\begin{align}
\u = \half \d_i \times \ps \d_i = \frac{1}{2L_0} \d_i \times \d_i'
= \frac{1}{L_0}\bar \u \quad \textnormal{and} \quad v_3 = \vert
\ps \r \vert = \vert \r' \vert = \bar v_3. \label{dimlessu}
\end{align}

Utilizing the time varying generic nodal displacements introduced
in (\ref{dynamicnda}) and (\ref{dynamicndb}) instead of the static
generic nodal displacements introduced in (\ref{staticnda}) and
(\ref{staticndb}) respectively, the time varying generic
displacements at any point within the element can be expressed as
nonlinear functions of the length parameter $\sigma$ and the nodal
displacement vector $\vq^e(\tau)$. Based on the nonlinear shape
functions derived in Section 4, we have
\begin{align} \left\{ \begin{array}{l}
\, \bar x = \hat x_1(\sigma, \tau) + \hat x_2(\sigma, \tau)
+ \hat x_3(\sigma, \tau),  \\
\, \bar y = \hat y_1(\sigma, \tau) + \hat y_2(\sigma, \tau)
+ \hat y_3(\sigma, \tau),  \\
\, \bar z = \sigma + \hat z_1(\sigma, \tau) + \hat z_2(\sigma,
\tau) + \hat z_3(\sigma, \tau),  \\
\, \varphi = \hat \varphi_1(\sigma, \tau) + \hat \varphi_2(\sigma,
\tau) + \hat \varphi_3(\sigma, \tau).
\end{array} \right. \label{genericd}
\end{align}
where the $i$th terms $\hat x_i$, $\hat y_i$, $\hat z_i$ and $\hat
\varphi_i$ are $i$th order functions of the nodal displacement
vector $\vq^e(\tau)$. For example, based on (\ref{curve1}) the 1st
order terms are
\begin{align}  \left\{ \begin{array}{ll}
\hat x_1(\sigma, \tau) =  X_a(\tau) & + \;\Phi_{ya}(\tau) \sigma -
(3X_a(\tau) - 3X_b(\tau) + 2l \Phi_{ya}(\tau) + l \Phi_{yb}(\tau))
\displaystyle  \frac{\sigma^2}{l^2} \\
& + \;(2X_a(\tau) - 2X_b(\tau) + l \Phi_{ya}(\tau) + l
\Phi_{yb}(\tau)) \displaystyle
\frac{\sigma^3}{l^3} \\
\hat y_1(\sigma, \tau) = Y_a(\tau) & - \;\Phi_{xa}(\tau) \sigma
-(3Y_a(\tau) - 3Y_b(\tau) - 2l
\Phi_{xa}(\tau) - l \Phi_{xb}(\tau)) \displaystyle  \frac{\sigma^2}{l^2} \\
& + \;(2Y_a(\tau) - 2Y_b(\tau) - l \Phi_{xa}(\tau) - l
\Phi_{xb}(\tau) )
\displaystyle  \frac{\sigma^3}{l^3} \\
\hat z_1(\sigma, \tau) =  Z_a(\tau) & + \;(Z_b(\tau) - Z_a(\tau))
\displaystyle  \frac{\sigma}{l}\\
\hat \phi_1(\sigma, \tau) = \Phi_{za}(\tau) & + \;(\Phi_{zb}(\tau)
- \Phi_{za}(\tau)) \displaystyle  \frac{\sigma}{l}. \end{array}
\right. \label{curve2}
\end{align}
The high order terms (up to third order), as indicated in Section
4, can be easily obtained using a MAPLE program \cite{LW}.
Consequently, the time varying generic displacements at any point
within the element can be written as
\begin{align}
\bar x = \bar x(\sigma, \vq^e(\tau)), \quad \bar y = \bar
y(\sigma, \vq^e(\tau)), \quad \bar z = \bar z(\sigma,
\vq^e(\tau)), \quad \varphi = \varphi(\sigma, \vq^e(\tau)).
\label{tvd}
\end{align}
This leads $\bar \r = \bar \r(\sigma, \vq^e(\tau))$. Moreover,
from (\ref{mu123}), we have
\begin{align}
\nu_1 = \frac{\bar
x'(\sigma,\vq^e(\tau))}{\abs{\bar\r'(\sigma,\vq^e(\tau))}} =
\nu_1(\sigma, \vq^e(\tau)), \quad \nu_2 = \frac{\bar
y'(\sigma,\vq^e(\tau))}{\abs{\bar\r'(\sigma,\vq^e(\tau))}} =
\nu_2(\sigma, \vq^e(\tau)) \label{tvnu}
\end{align}
Substituting $\nu_1(\sigma,\vq^e(\tau))$,
$\nu_2(\sigma,\vq^e(\tau))$ and $\varphi(\sigma,\vq^e(\tau))$ into
the expressions (\ref{di1})--(\ref{di3}) yields
\begin{align}
\d_i = \d_i(\sigma,\vq^e(\tau)), \hspace{1cm} i=1,2,3. \label{diq}
\end{align}
Similarly, from (\ref{rota}), we have
\begin{align}
\phi_x = \phi_x(\sigma,\vq^e(\tau)), \; \quad \phi_y =
\phi_y(\sigma,\vq^e(\tau)), \; \quad \phi_z =
\phi_z(\sigma,\vq^e(\tau)).
\end{align}
It follows from (\ref{dimlessw}), (\ref{dimlessu}) and (\ref{diq})
that
\begin{align}
\bar \w = \frac{1}{2} \d_i \times \dot \d_i = \bar
\w(\sigma,\vq^e(\tau)), \quad \bar \u = \frac{1}{2} \d_i \times
\d_i' = \bar \u(\sigma,\vq^e(\tau))
\end{align}
Therefore, the kinetic energy density (\ref{k1}) and the potential
energy density (\ref{p1}) are expressed as
\begin{align}
\T = \T(\sigma,\vq^e(\tau), \dot \vq^e(\tau)), \quad \U =
\U(\sigma,\vq^e(\tau)).
\end{align}
Then, the Lagrangian defined in the classical form $\mathscr{L} =
\mathscr{T} - \mathscr{V}$ are obtained as
\begin{align}
\mathscr{L}(\vq^e, \dot \vq^e) = \mathscr{T}(\vq^e, \dot \vq^e)  -
\mathscr{V}(\vq^e)  = \int_0^l (\,\T(\sigma,\vq^e, \dot \vq^e)-
\U(\sigma,\vq^e)\,) \,L_0\,d\sigma. \label{Lagrangian}
\end{align}

So far we have not precisely defined the type of loading. Let us
assume that a load acting on the element is composed from three
additive parts. The first one is the interaction of the neighbored
elements. The second one is the external point (concentrated)
loadings acting on the nodes. The last one represents a
distributed load with fixed direction and prescribed intensity as
mentioned in Section 2. In keeping with the load definitions in
the principle of virtual work, the total load has to be defined
with respect to the inertial basis because the generalized nodal
displacements are defined with respect to them. Thus, let us
denote
\begin{align}
\vf^i_a(\tau) = \left[ \begin{array}{c} f_{xa}^i(\tau) \\ f_{ya}^i(\tau) \\
f_{za}^i(\tau) \end{array} \right], \quad \vf^i_b = \left[
\begin{array}{c} f_{xb}^i(\tau) \\ f_{yb}^i(\tau) \\
f_{zb}^i(\tau) \end{array} \right],
 \quad \quad \vl^i_a = \left[
\begin{array}{c} l_{xa}^i(\tau) \\ l_{ya}^i(\tau) \\
l_{za}^i(\tau) \end{array} \right],
\quad \vl^i_b = \left[
\begin{array}{c} l_{xb}^i(\tau) \\ l_{yb}^i(\tau) \\
l_{zb}^i(\tau) \end{array} \right]
\end{align}
be the interaction force and torque vector at nodes $\sigma = 0$
and $\sigma = l$ respectively.

Similarly, the external point loadings are expressed as
\begin{align}
\vf^c_a(\tau) = \left[ \begin{array}{c} f_{xa}^c(\tau) \\ f_{ya}^c(\tau) \\
f_{za}^c(\tau) \end{array} \right], \quad \vf^c_b(\tau) = \left[
\begin{array}{c} f_{xb}^c(\tau) \\ f_{yb}^c(\tau) \\
f_{zb}^c(\tau) \end{array} \right],
 \quad \quad \vl^c_a = \left[
\begin{array}{c} l_{xa}^c(\tau) \\ l_{ya}^c(\tau) \\
l_{za}^c(\tau) \end{array} \right],
\quad \vl^c_b = \left[
\begin{array}{c} l_{xb}^c(\tau) \\ l_{yb}^c(\tau) \\
l_{zb}^c(\tau) \end{array}
\right]\;,
\end{align}
while the distributed forces and torques may be expressed as
\begin{align}
\boldsymbol{\xi}^d = \left[ \begin{array}{c} \xi_x^d(\sigma, \tau) \\
\xi_y^d(\sigma, \tau) \\
\xi_z^d(\sigma, \tau) \end{array} \right], \quad
\boldsymbol{\eta}^d =
\left[ \begin{array}{c} \eta_x^d(\sigma, \tau) \\ \eta_y^d(\sigma, \tau) \\
\eta_z^d(\sigma, \tau) \end{array} \right].
\end{align}
The virtual work done by the distributed load has the form
\begin{align}
\delta \mathscr{W}^d & = \int_0^l \left( \xi_x^d \delta \bar x +
\xi_y^d \delta \bar y + \xi_z^d \delta \bar z + \eta_x^d \delta
\phi_x + \eta_y^d \delta \phi_y + \eta_z^d \delta \phi_z
\right)\,L_0\, d\sigma \nonumber \\[12pt]
& = \int_0^l \left( \xi_x^d \frac{\partial \bar x(\sigma,
\vq^e)}{\partial \vq^e} + \xi_y^d \frac{\partial \bar y(\sigma,
\vq^e)}{\partial \vq^e} + \xi_z^d \frac{\partial \bar z(\sigma,
\vq^e)}{\partial \vq^e} + \right. \nonumber \\[12pt]
& \left. \;\;\;\; + \eta_x^d \frac{\partial \phi_x(\sigma,
\vq^e)}{\partial \vq^e} + \eta_y^d \frac{\partial \phi_y(\sigma,
\vq^e)}{\partial \vq^e} + \eta_z^d \frac{\partial \phi_z(\sigma,
\vq^e)}{\partial \vq^e} \right)\, \delta \vq^e\, L_0\, d\sigma
\end{align}
For the sake of convenience, let
\begin{align}
& \vf^{ie}(\tau) = \left[ \begin{array}{c} \vf_a^i(\tau) \\ \vl_a^i(\tau) \\
\vf_b^i(\tau) \\ \vl_b^i(\tau) \end{array} \right], \quad
\vf^{ce}(\tau) = \left[
\begin{array}{c}
\vf_a^c(\tau) \\ \vl_a^c(\tau) \\
\vf_b^c(\tau) \\ \vl_b^c(\tau) \end{array} \right], \\[12pt]
& \vf^{de}(\tau, \vq^e) = \int_0^l \left( \xi_x^d(\tau)
\frac{\partial \bar x(\sigma, \vq^e)}{\partial \vq^e} + \cdots +
\eta_z^d(\tau) \frac{\partial \phi_z(\sigma, \vq^e)}{\partial
\vq^e} \right)^T\, L_0\, d\sigma.
\end{align}
Then, the total virtual work done by the three additive parts are
\begin{align}
\delta \mathscr{W} = (\vf^{ie} + \vf^{ce} + \vf^{de})^T \cdot
\delta \vq^e. \label{workdone}
\end{align}

Substituting (\ref{Lagrangian}) and (\ref{workdone}) into
(\ref{Hamilton}), taking variations using the chain rule, and
integrating by parts, yield the generalized Lagrange equations of
motion for the Cosserat rod element:
\begin{align}
\frac{d}{d\tau} \left( \frac{\partial L}{\partial \dot q_j}
\right) - \frac{\partial L}{\partial q_j} = f_j^{ie}(\tau) +
f_j^{ce}(\tau) + f_j^{de}(\tau, \vq^e)
\end{align}

For a general configuration with nonzero generic nodal
displacements $\vq^e$, the ordinary differential equations of
motion with up to third order nonlinearities of displacements and
first order kinetic terms can be obtained as
\begin{align}
\vM^e \ddot \vq^e + \vK^e \vq^e + \vg^e(\vq^e) = \vf^{ie}(\tau) +
\vf^{ce}(\tau) + \vf^{de}(\tau, \vq^e), \label{odecre}
\end{align}
where $\vM^e$ and $\vK^e$ are mass and (linear) stiffness matrices
of the element model, $\vg^e(\vq^e)$ is a nonlinear vector with
quadratic and cubic terms of $\vq^e$. Since the mass of a typical
rod, such as the springlike support component of MEMS, is very
small comparing with the mass of the main device in practice only
the first order kinetic terms are reserved in Equation
(\ref{odecre}).

The detailed expressions of $\vM^e$, $\vK^e$ and $\vg^e(\vq^e)$
have been implemented in a MAPLE program \cite{LW}. For the sake
of illustration, the explicit expressions of $\vM^e$, $\vK^e$ and
$\vg^e(\vq^e)$ for a cantilever beam as a special Cosserat rod
element are listed in Appendix.

\section{Dynamical Responses of rods by Cosserat Rod Elements}

\subsection{Assembly of equations of motion for whole system}

We could analyze all of the types of systems consist of a set of
interconnected components described in the introduction by using
Cosserat rod elements for the deformable parts or subdivided
members. Two- and three-dimensional frame structures require
rotation-of-axes transformation for actions and displacements. For
the sake of convenience, in this section we shall examine only the
type of structure which is aligned with reference axes, using
properties of the Cosserat rod element developed in the preceding
sections. The analysis of the response of a number of complex
structures is beyond the scope of this paper and will be presented
in future publications.

After stiffness, mass, and actual or equivalent nodal loads for
individual Cosserat rod element are generated, we can assembly
them to form the equations of motion for a whole system. We define
global displacement vector $\vq$ holding the displacement
variables for all mesh nodes, such that
\begin{align}
\vq = \left[ X_1\; Y_1\; Z_1\; \Phi_{x1}\; \Phi_{y1}\; \Phi_{z1}\;
X_2\; Y_2\; Z_2\; \Phi_{x2}\; \Phi_{y2}\; \Phi_{z2}\; \cdots
\right]^T.
\end{align}

The equations of motion for the whole system can be constructed by
simply adding the contributions from all the elements. In this
way, expanding the matrix or operator for each individual element
to make them the same size as the system matrices or operators, we
have
\begin{align}
\vM = \sum_{i=1}^{n_e} \vM_i^e, \qquad \vK = \sum_{i=1}^{n_e}
\vK^e_i, \label{sys1}
\end{align}
and
\begin{align}
\vg(\vq) = \sum_{i=1}^{n_e} \vg^e_i(\vq), \qquad \vf^c(\tau) =
\sum_{i=1}^{n_e} \vf_i^{ce}(\tau), \qquad \vf^d(\tau,\vq) =
\sum_{i=1}^{n_e} \vf_i^{de}(\tau, \vq). \label{sys2}
\end{align}
where $n_e$ is the number of elements. In Equation (\ref{sys1})
$\vM$ and $\vK$ represent the system mass matrix and the system
(linear) stiffness matrix. Similarly, the action vectors
$\vf^c(\tau)$ and $\vf^d(\tau, \vq)$ are actual and equivalent
nodal loads for the whole system. The contributions from the
interaction forces and torques from all the elements must be of
balance and the total action must be vanished. Then the undamped
equations of motion for the assembled system become
\begin{align}
\vM \ddot \vq + \vK \vq + \vg(\vq) = \vf^c(\tau) +
\vf^d(\tau,\vq). \label{sys3}
\end{align}
This equation gives the system equations of motion all nodal
displacements, regardless of whether they are free or restricted.

In preparation for solving the nonlinear dynamic equations
(\ref{sys3}), as in the standard finite element procedure, we
rearrange and partition it as follows
\begin{align}
\left[ \begin{array}{cc} \vM_{ff}\; & \;\vM_{fr} \\ \vM_{rf} &
\vM_{rr} \end{array} \right] \; \left[  \begin{array}{c} \;\ddot
\vq_f \\ \ddot \vq_r \end{array} \right] + & \left[
\begin{array}{cc} \vK_{ff}\; & \;\vK_{fr} \\ \vK_{rf} &
\vK_{rr} \end{array} \right] \; \left[  \begin{array}{c} \;\vq_f
\\ \vq_r \end{array} \right]  + \left[
\begin{array}{c} \;\vg_f(\vq_f, \vq_r) \\ \vg_r(\vq_f, \vq_r) \end{array} \right]
\nonumber \\[10pt]
= & \left[ \begin{array}{c} \;\vf_f^c(\tau) + \vf_f^d(\tau,\vq_f,\vq_r) \\
\vf_r^c(\tau) + \vf_r^d(\tau,\vq_f,\vq_r) \end{array} \right],
\label{sys4}
\end{align}
in which the subscript $f$ refers to free nodal displacements
while the subscript $r$ denotes restrained nodal displacements. If
the support motions (at constraints) are zero, the equation
(\ref{sys4}) can be simplified to
\begin{align}
\vM_{ff} \ddot \vq_f + \vK_{ff} \vq_f + \vg_f(\vq_f) =
\vf_f^c(\tau) + \vf_f^d(\tau,\vq_f), \label{sys5}
\end{align}
and
\begin{align}
\vM_{rf} \ddot \vq_f + \vK_{rf} \vq_f + \vg_r(\vq_f) =
\vf_r^c(\tau) + \vf_r^d(\tau,\vq_f), \label{sys6}
\end{align}
which can be used for solving the free displacements $\vq_f(\tau)$
and support actions $\vf_r^c(\tau)$, respectively.

\subsection{Simulation results and discussion for a simple cantilever}

\begin{figure}[h!]
\centerline{
\includegraphics[height=5.0cm]{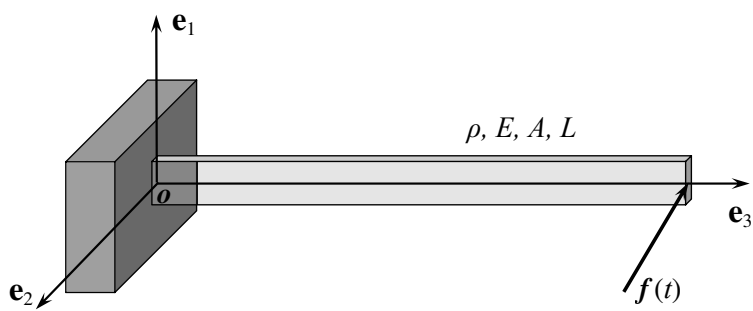}}
\caption{Schematic of a simple Cantilever} \label{figCmode2}
\end{figure}

A cantilever, as shown in Figure \ref{figCmode2}, is now presented
as a simple example to demonstrate high accuracy and excellent
performance of the proposed Cosserat rod elements. Numerical
calculations based on (\ref{sys5}) are carried out for a uniform
horizontal cantilever of length $L = 0.3$m, of constant cross
section with width $B=0.01$m and thickness $D=0.005$m. The mass
density and the Young's modulus are assumed to be $\rho = 3.0
\times 10^3$kg/m$^3$ and $E = 2.08 \times 10^8$Pa.

Dividing the cantilever into $n_e$ elements of equal length, we
can establish the nonlinear differential equations of motion
(\ref{sys5}) for solving the free displacements. In what follows,
the natural frequencies of the linearized system are studied and
used to compare with those derived from the classical beam theory
presented in textbooks (see, for example \cite{Ginsberg}), and
numerical simulations for the responses of the nonlinear dynamical
system (\ref{sys5}) under external harmonic excitations are
performed with Matlab.
\begin{center}
\begin{tabular}{p{1.6cm}|p{5.94cm}|p{5.94cm}}
\multicolumn{3}{l} {{\bf Table 1}\;\; Flexural natural frequencies
based on CRE approach and exact Continua method} \\[4pt]
\hline \hspace{0.5cm}$\omega_i$ & Flexural frequencies in
$\e_1$-$\e_3$ plane
& Flexural frequencies in $\e_2$-$\e_3$ plane\\[-1pt]
\end{tabular}
\begin{tabular}{p{1.6cm}|p{1.7cm} p{1.7cm} p{1.7cm}|p{1.7cm} p{1.7cm} p{1.7cm} }
 (rad/sec)  &  CRE &  CBT  &
$\vert$Error$\vert$ ($\%$) &  CRE  &  CBT  & $\vert$Error$\vert$ ($\%$) \\
\hline \hspace{0.5cm}1 & 29.7607 & 29.7665 & 0.0197 & 14.8827 & 14.8833 & 0.0036 \\
\hspace{0.5cm}2 & 186.358 & 186.544 & 0.0995 & 93.2838 & 93.2718 & 0.0129 \\
\hspace{0.5cm}3 & 522.329 & 522.329 & 0.0000 & 261.868 & 261.164 & 0.2692\\
\hspace{0.5cm}4 & 1028.68 & 1023.56 & 0.5005 & 516.914 & 511.778 & 1.0035\\
\hspace{0.5cm}5 & 1707.74 & 1692.01 & 0.5155 & 857.104 & 846.007 & 1.3118\\
 \hline
\end{tabular}
\end{center}
First, the flexural natural frequencies calculated in terms of the
linearized equations of the nonlinear system (\ref{sys5}) obtained
by Cosserat element approach, together with the theoretical
results obtained by employing the classical beam theory (CBT) are
given in Table 1. The flexural natural frequencies in both
$\e_1$-$\e_3$ plane and $\e_2$-$\e_3$ plane, based on the CRE
approach, showed their excellent convergency (the corresponding
results listed in Table 1 are found, when only five Cosserat rod
elements are used).

\begin{figure}[h!]
\centerline{
\includegraphics[height=6.0cm]{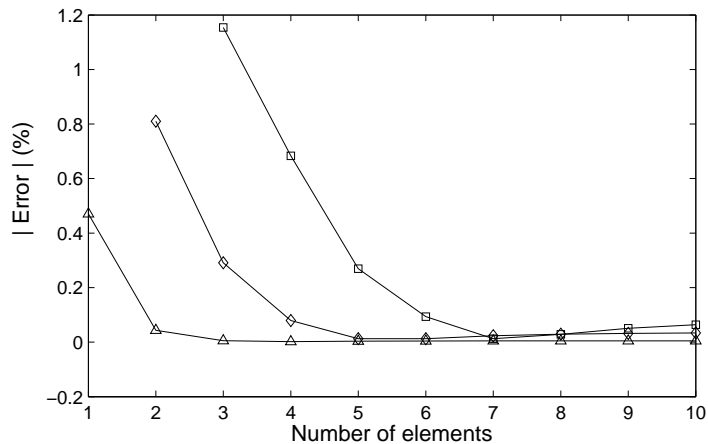}}
\caption{Convergency test for the first flexural frequencies in
$\e_2$-$\e_3$ plane. --$\vartriangle$--, $\omega_1$;
--$\lozenge$--, $\omega_2$; --{\scriptsize $\square$}--,
$\omega_3$.} \label{errorx}
\end{figure}

Figure \ref{errorx} represents the CRE convergency tests
corresponding to the first three flexural natural frequencies in
$\e_2$-$\e_3$ plane of the rod. As can be seen for the first
frequency, the $\vert$error$\vert$ is found to be very small ($\le
0.1\%$) even when only two elements are used. In fact the
$\vert$error$\vert$ for the first frequency is only $0.4535\%$
when just one element is used. For the second and third natural
frequencies, the results are converging with approximately $0.1\%$
error, when six elements are used.

In the second part of this example, based on the derived nonlinear
system (\ref{sys5}), numerical simulations are performed to
investigate the dynamic responses of the cantilever under harmonic
excitations. The differential equations of motion are full coupled
by the nonlinear terms and could exhibit internal resonance
introduced by the nonlinearities. They also exhibit external
resonances when the external excitation is periodic and the
frequency of a component of its Fourier series is near one of the
natural frequencies of the system, or near a multiple of a natural
frequencies. The detailed analysis of complex dynamic behavior,
such as bifurcation and chaos, of the system is not the main focus
of this paper. We only compare here the responses of the system,
when different number of elements are used.

\begin{figure}[h!]
\centerline{
\includegraphics[height=8.5cm]{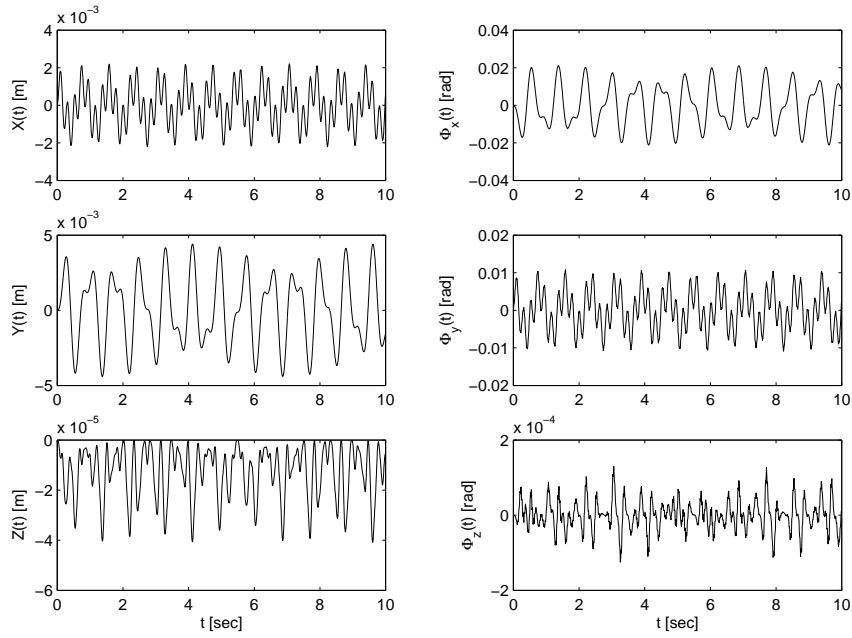}}
\caption{Displacement time histories of the rod with external
loads $f_x(t) = 0.01 \cos(8*t)$, $f_y(t) = 0.005 \sin(8*t)$ and
zero initial conditions: two elements case.} \label{twoe}
\end{figure}
\begin{figure}[h!]
\centerline{
\includegraphics[height=8.5cm]{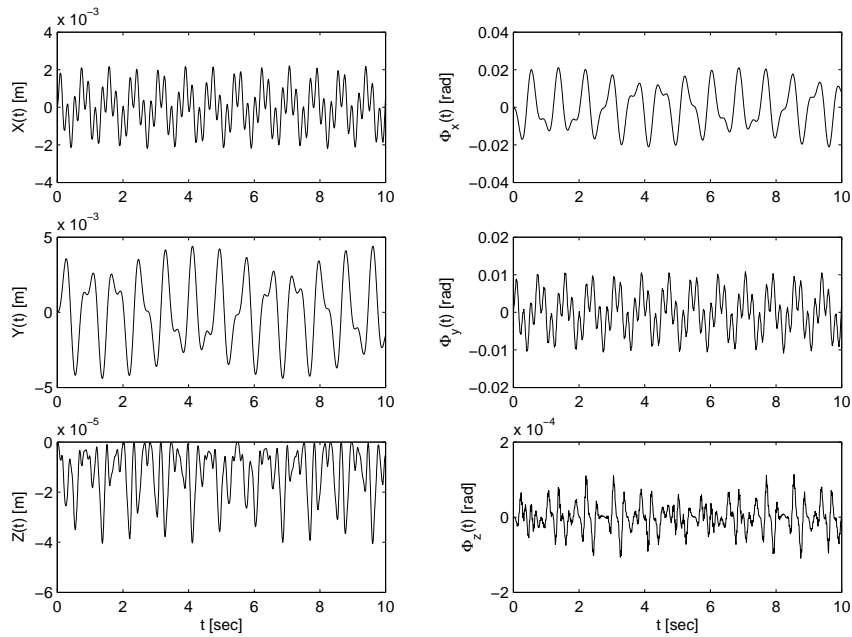}}
\caption{Displacement time histories of the rod with external
loads $f_x(t) = 0.01 \cos(8*t)$, $f_y(t) = 0.005 \sin(8*t)$ and
zero initial conditions: ten elements case.} \label{tene}
\end{figure}

\begin{figure}[h!]
\centerline{
\includegraphics[height=11.50cm]{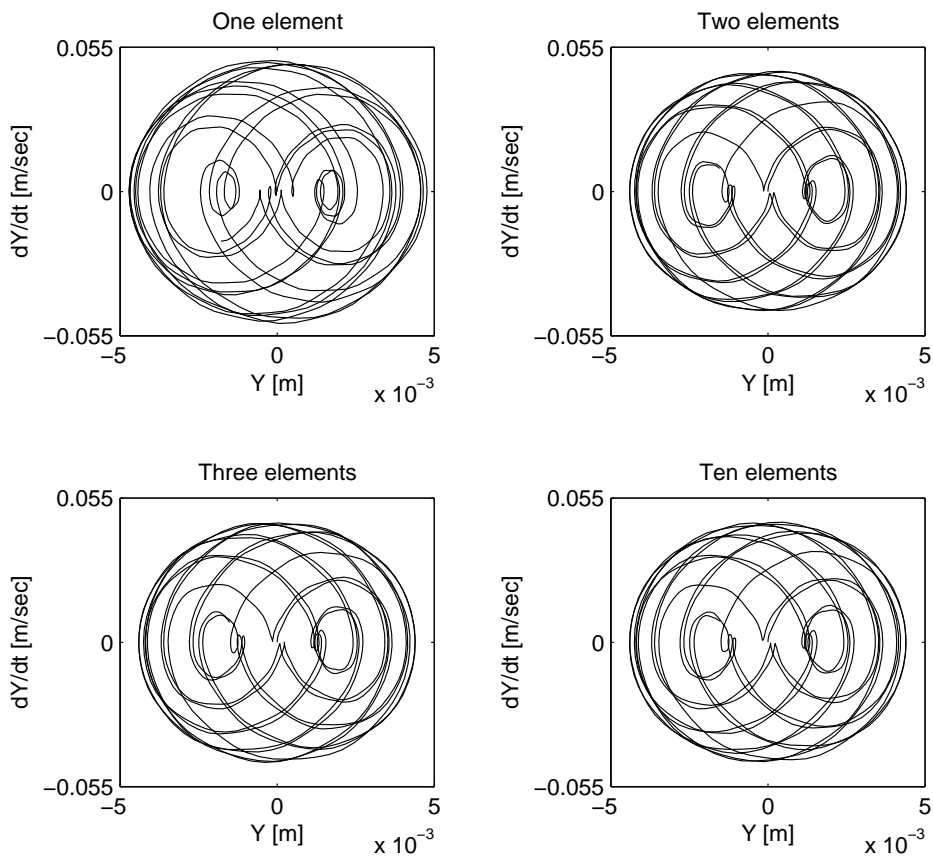}}
\caption{Phase plane diagram of $Y$--$\dot Y$ with external loads
$f_x(t) = 0.01 \cos(8*t)$, $f_y(t) = 0.005 \sin(8*t)$ and zero
initial conditions.} \label{ottte}
\end{figure}

The displacement and angular time histories of the free end of the
cantilever under external loads $f^c_x(t) = 0.01 \cos(8t)$,
$f^c_y(t) = 0.005 \sin(8t)$ and at zero initial conditions are
shown in Figure \ref{twoe} when two elements are used and in
Figure \ref{tene} when ten elements are used, respectively. It is
interesting to note that amplitudes and periods of the responses
are very closed in this two situations. To enhance this
observation, the phase plane diagrams for $Y(t)$--$\dot Y(t)$ in
four different cases, namely one element, two elements, three
elements and ten elements, are plotted in Figure \ref{ottte} (a),
(b), (c) and (d), respectively. Comparing the four diagrams in
Figure \ref{ottte} shows that the modal when two or three elements
are used can exhibit almost the same behavior of the model when
ten elements are used.

According to the analysis of natural frequencies and the analysis
of harmonic responses of the established nonlinear dynamic
systems, we believe, in practical engineering problem, especial
for the structure composed of springlike flexural components such
as the device in MEMS, only a few Cosserat rod elements are needed
to model a flexural component. For the very slender flexural
components, we can even use only one element to model such a
component.

\section{Conclusion}

A Cosserat rod element formulation for the modelling of
3-dimensional dynamics of slender structures has been proposed in
this paper. The modelling strategy of this new approach employed
the exact nonlinear kinematic relationships in the sense of
Cosserat theory, and adopted the Bernoulli hypothesis. Finite
displacements and rotations as well as finite extensional,
torsional, and bending strains are accounted for. The Kirchoff
constitutive relations, which provide an adequate description of
elastic properties in terms of a few elastic moduli, are adopted.
A deformed configuration of the rod is described by the
displacement vector of the deformed centroid curves and an
orthonormal moving frame, rigidly attached to the cross-section of
the rod. The position of the moving frame relative to the inertial
frame is specified by the rotation matrix, parametrized by a
rotational vector. The approximation solutions of the nonlinear
partial differential equations of motion in quasi-static sense are
chosen as the shape functions with up to third order nonlinear
terms of generic nodal displacements. This lends the approach very
well to achieve higher accuracy of the dynamic responses of the
model by dividing the slender rod into a few elements. Based on
the Lagrangian constructed by the Cosserat kinetic energy and
strain energy expressions, the principle of virtual work is
employed to derive the ordinary differential equations of motion
with third order nonlinear generic nodal displacements.

A cantilever as a simple example has been presented to illustrate
the use of the formulation developed here to obtain the lower
order nonlinear ordinary differential equations of motion of a
given structure. The natural frequency analysis for the linearized
equations and the numerical simulation analysis for the nonlinear
model show that in practical engineering problem, especial for the
structure composed of springlike flexural components such as the
device in MEMS, only a few Cosserat rod elements are needed to
model a flexural component.

The mathematical simplicity when formulating deformable components
enables more convenient for modelling the multibody systems that
consist of interconnected rigid and deformable components. The
Cosserat rod element approach therefore is feasible to be used to
capture the most significant characteristics of a multi- rigid and
deformable body system in a few variables governed by nonlinear
ordinary differential equations of motion.

As the first step to present the Cosserat rod element approach, we
have limited our attention to the modelling of Cosserat rod
elements in which the effect of shear has been neglected. The
extension of the present formulation to the modelling of more
general Cosserat rod elements in which the finite extensional,
torsional, bending strains as well as shear are accounted for is
highly desirable.

\vspace{0.5cm} \noindent {\bf Acknowledgements} The authors are
grateful to the EPSRC (Computational Engineering Mathematics
Programme) and the EC(Framework Programme) for financial support
in this study.

\newpage

\noindent {\bf \Large Appendix}

Let us assume that a uniform cantilever beam of length $L$, of
constant cross section with area $A$ and density $\rho$, is fixed
at $s=0$ and free at $s=L$. In this case, we have $\vq_a = 0$,
thus $\vq^e = \vq_b$. Consequently, $\vM^e, \vK^e$ become $6\times
6$ matrices, and $\vg^e(\vq^e)$ is a six dimensional nonlinear
vectorial functions of $\vq^e = \vq_b$. They are
\begin{align}
\vM^e = \left[  \begin{array}{cccccc} \displaystyle \frac{13\mu
l^2 + 42 I_{22}}{35l} & 0 & 0 & 0 & \displaystyle -
\frac{I_{22}}{10} - \frac{11\mu l^2}{210} & 0 \\
0 & \displaystyle \frac{13\mu l^2 + 42 I_{11}}{35l} & 0 &
\displaystyle \frac{I_{11}}{10} + \frac{11\mu l^2}{210} & 0 & 0\\
0 & 0 & \; \displaystyle \frac{\mu l}{3} \; & 0 & 0 & 0\\
0 & \displaystyle \frac{I_{11}}{10} + \frac{11\mu l^2}{210} & 0 &
\displaystyle \frac{2I_{11}l}{15} + \frac{\mu l^3}{105} & 0 & 0\\
\displaystyle - \frac{I_{22}}{10} - \frac{11\mu l^2}{210} & 0 & 0
& 0 & \displaystyle
\frac{2I_{22}l}{15} + \frac{\mu l^3}{105} & 0\\
0 & 0 & 0 & 0 & 0 & \; \displaystyle \frac{I_{33}l}{3} \;
\end{array} \right]\;,
\end{align}
\begin{align}
\vK_c^e = \left[  \begin{array}{cccccc} \displaystyle
\frac{12J_{22}}{l^3} & 0 & 0 & 0 & \displaystyle -
\frac{6J_{22}}{l^2} & 0 \\
0 & \displaystyle \frac{12J_{11}}{l^3} & 0 &
\displaystyle \frac{6J_{11}}{l^2} & 0 & 0\\
0 & 0 & \; \displaystyle \frac{K_{33}}{l} \; & 0 & 0 & 0\\
0 & \displaystyle \frac{6J_{11}}{l^2} & 0 &
\displaystyle \frac{4J_{11}}{l} & 0 & 0\\
\displaystyle - \frac{6J_{22}}{l^2} & 0 & 0 & 0 & \displaystyle
\frac{4J_{22}}{l} & 0\\
0 & 0 & 0 & 0 & 0 & \; \displaystyle \frac{J_{33}}{l} \;
\end{array} \right]\;,
\end{align}
and
\begin{align}
g_1(\vq^e)\, = \; & g_{1,1} X_bZ_b + g_{1,2} Y_b \Phi_{zb} +
g_{1,3}Z_b \Phi_{yb} + g_{1,4} \Phi_{xb} \Phi_{zb} + g_{1,5} X_b^3
+ g_{1,6} X_b^2 \Phi_{yb} \nonumber\\
&  + g_{1,7} X_b Y_b^2 + g_{1,8}X_b Y_b \Phi_{xb} + g_{1,9}X_b
Z_b^2 + g_{1,10} X_b \Phi_{xb}^2 + g_{1,11} X_b \Phi_{yb}^2  \nonumber\\
& + g_{1,12} X_b \Phi_{zb}^2 + g_{1,13} Y_b^2 \Phi_{yb} + g_{1,14}
Y_b Z_b \Phi_{zb} + g_{1,15}Y_b \Phi_{xb} \Phi_{yb} + g_{1,16}
Z_b^2 \Phi_{yb}  \nonumber\\
& + g_{1,17} Z_b \Phi_{xb} \Phi_{zb} + g_{1,18} \Phi_{xb}^2
\Phi_{zb} + g_{1,19} \Phi_{yb}^3 + g_{1,20} \Phi_{yb}
\Phi_{zb}^2\;,
\\[5pt]
g_2(\vq^e)\, = \; & g_{2,1} X_b \Phi_{zb} + g_{2,2} Y_b Z_b +
g_{2,3} Z_b \Phi_{xb} + g_{2,4} \Phi_{yb} \Phi_{zb} + g_{2,5}
X_b^2 Y_b + g_{2,6} X_b^2 \Phi_{xb}, \nonumber\\
& + g_{2,7} X_b Y_b \Phi_{yb} + g_{2,8} X_b Z_b \Phi_{zb} +
g_{2,9} X_b \Phi_{xb} \Phi_{yb} + g_{2,10}Y_b^3 + g_{2,11}Y_b^2
\Phi_{xb} \nonumber\\
& + g_{2,12} Y_b Z_b^2 + g_{2,13}Y_b \Phi_{xb}^2 + g_{2,14}Y_b
\Phi_{yb}^2 + g_{2,15}Y_b \Phi_{zb}^2 + g_{2,16}Z_b^2 \Phi_{xb} \nonumber\\
& + g_{2,17} Z_b \Phi_{yb} \Phi_{zb} + g_{2,18} \Phi_{xb}^3 +
g_{2,19} \Phi_{xb} \Phi_{yb}^2 + g_{2,20} \Phi_{xb} \Phi_{zb}^2\;,
\end{align}
\begin{align}
g_3(\vq^e)\, = \; & g_{3,1} X_b^2 + g_{3,2} X_b \Phi_{yb} +
g_{3,3} Y_b^2 + g_{3,4} Y_b \Phi_{xb} + g_{3,5} \Phi_{xb}^2 +
g_{3,6}\Phi_{yb}^2
 + g_{3,7} X_b^2 Z_b \nonumber\\
& + g_{3,8} X_b Y_b \Phi_{zb} + g_{3,9} X_b Z_b \Phi_{yb} +
g_{3,10} X_b \Phi_{xb} \Phi_{zb} + g_{3,11}Y_b^2
Z_b + g_{3,12}Y_b Z_b \Phi_{xb} \nonumber\\
& + g_{3,13} Y_b \Phi_{yb} \Phi_{zb} + g_{3,14} Z_b \Phi_{xb}^2 +
g_{3,15}Z_b \Phi_{yb}^2 + g_{3,16} \Phi_{xb} \Phi_{yb}
\Phi_{zb}\;,\\[5pt]
g_4(\vq^e)\, = \; & g_{4,1} X_b \Phi_{zb} + g_{4,2} Y_b Z_b +
g_{4,3} Z_b \Phi_{xb} + g_{4,4} \Phi_{yb} \Phi_{zb} + g_{4,5}
X_b^2 Y_b + g_{4,6} X_b^2 \Phi_{xb} \nonumber\\
& + g_{4,7} X_b Y_b \Phi_{yb} + g_{4,8} X_b Z_b \Phi_{zb} +
g_{4,9} X_b \Phi_{xb} \Phi_{yb} + g_{4,10}Y_b^3 + g_{4,11}Y_b^2
\Phi_{xb} \nonumber\\
& + g_{4,12} Y_b Z_b^2 + g_{4,13}Y_b \Phi_{xb}^2 + g_{4,14}Y_b
\Phi_{yb}^2 + g_{4,15}Y_b \Phi_{zb}^2 + g_{4,16}Z_b^2 \Phi_{xb} \nonumber\\
& + g_{4,17} Z_b \Phi_{yb} \Phi_{zb} + g_{4,18} \Phi_{xb}^3 +
g_{4,19} \Phi_{xb} \Phi_{yb}^2 + g_{4,20} \Phi_{xb} \Phi_{zb}^2\;,
\\[5pt]
g_5(\vq^e)\, = \; & g_{1,1} X_bZ_b + g_{5,2} Y_b \Phi_{zb} +
g_{5,3}Z_b \Phi_{yb} + g_{5,4} \Phi_{xb} \Phi_{zb} + g_{5,5} X_b^3
+ g_{5,6} X_b^2 \Phi_{yb} \nonumber\\
&  + g_{5,7} X_b Y_b^2 + g_{5,8}X_b Y_b \Phi_{xb} + g_{5,9}X_b
Z_b^2 + g_{1,10} X_b \Phi_{xb}^2 + g_{1,11} X_b \Phi_{yb}^2  \nonumber\\
& + g_{5,12} X_b \Phi_{zb}^2 + g{5,13} Y_b^2 \Phi_{yb} + g_{5,14}
Y_b Z_b \Phi_{zb} + g_{5,15}Y_b \Phi_{xb} \Phi_{yb} + g_{5,16}
Z_b^2 \Phi_{yb}  \nonumber\\
& + g_{5,17} Z_b \Phi_{xb} \Phi_{zb} + g_{5,18} \Phi_{xb}^2
\Phi_{zb} + g_{5,19} \Phi_{yb}^3 + g_{5,20} \Phi_{yb} \Phi_{zb}^2\;, \\[5pt]
g_6(\vq^e)\, = \; & g_{6,1} X_b Y_b + g_{6,2} X_b \Phi_{xb} +
g_{6,3} Y_b \Phi_{yb} + g_{6,4} \Phi_{xb} \Phi_{yb} + g_{6,5}
x_2^2 \Phi_{zb} + g_{6,6} X_b Y_b Z_b \nonumber\\
&  + g_{6,7} X_b Z_b \Phi_{xb} + g_{6,8} X_b \Phi_{yb} \Phi_{zb} +
g_{6,9} Y_b^2 \Phi_{zb} + g_{6,10} Y_b Z_b \Phi_{yb} + g_{6,11}
Y_b \Phi{xb} \Phi{zb} \nonumber\\
& + g_{6,12} Z_b \Phi_{xb} \Phi_{yb} +  g_{6,13} \Phi_{xb}^2
\Phi_{zb} + g_{6,14} \Phi_{yb}^2 \Phi_{zb}\;,
\end{align}
where, the coefficients of second order nonlinear terms are
specified as:
\begin{align} \begin{array}{ll}
g_{1,1} = 2 g_{3,1} = \displaystyle \frac{6(K_{33}l^2 -
20J_{22})}{5l^4},\hspace{2cm} & g_{1,2} = g_{2,1} = g_{6,1} =
\displaystyle \frac{6(J_{22} - J_{11})}{l^3},
\\[8pt]
g_{1,3} = g_{3,2} = - g_{5,1} = \displaystyle \frac{K_{33}l^2 -
60J_{22}}{10l^3} ,\qquad & g_{1,4} = g_{4,1} = - g_{6,2} =
\displaystyle \frac{4J_{11} - J_{22} - J_{33}}{l^2}, \\[8pt]
g_{2,2} = 2 g_{3,3} = \displaystyle \frac{6(K_{33}l^2 -
20J_{11})}{5l^4}, & g_{2,3} = g_{3,4} = g_{4,2} = \displaystyle
\frac{K_{33}l^2 - 60J_{11}}{10l^3}, \\[8pt]
g_{2,4} = - g_{5,2} = g_{6,3} = \displaystyle
\frac{J_{11}-4J_{22}+J_{33}}{l^2}, & g_{3,5} = g_{3,6} =
\frac{1}{2} g_{4,3} = - \frac{1}{2} g_{5,3} = \displaystyle
\frac{K_{33}}{15}, \\[8pt]
g_{4,4} = - g_{5,4} = g_{6,4} = \displaystyle
\frac{J_{11}-J_{22}}{l}. &
\end{array} \nonumber \end{align}
The coefficients of third order nonlinear terms are specified as :
\begin{align} \begin{array}{l}
g_{1,5} = \displaystyle \frac{ 18(7K_{33}^2l^4 -
160J_{22}K_{33}l^2 - 560J_{22}^2)}{175K_{33}l^7}\;, \\[8pt]
g_{1,6} = \displaystyle \frac{ 9(7K_{33}^2l^4 -
260J_{22}K_{33}l^2 - 3360J_{22}^2)}{350K_{33}l^6}\;, \\[8pt]
g_{1,7} = \displaystyle \frac{18(7K_{33}l^2 - 80(J_{11} +
J_{22}))}{175 l^5} - \frac{18(10K_{33}l^2(J_{11}- J_{22})^2 + 112
J_{11}J_{22}J_{33})}{35J_{33}K_{33}l^7}\;, \\[8pt]
g_{1,8} = \displaystyle \frac{ 3(7K_{33}l^2 - 480J_{11} +
220J_{22})}{175 l^4} - \frac{ 18(10K_{33}l^2(J_{11}-J_{22})^2 +
112J_{11}J_{22}J_{33})}{35J_{33}K_{33}l^6}\;,
\end{array} \nonumber
\end{align}

\begin{align} \begin{array}{l}
g_{1,9} = - \displaystyle \frac{ K_{33}^2l^4 + 840J_{22}K_{33} -
25200J_{22}^2}{700J_{22}l^5}\;,  \\[10pt]
g_{1,10} = \displaystyle \frac{14K_{33}l^2 -500J_{11} -80J_{22} +
175J_{33}}{175l^3} - \frac{52K_{33}l^2(J_{11}-J_{22})^2 +
504J_{11}J_{22}J_{33}}{35J_{33}K_{33}l^5}\;,\\[10pt]
g_{1,11} = \displaystyle \frac{ 63K_{33}^2l^4 - 520J_{22}K_{33}l^2
- 38640J_{22}^2}{700K_{33}l^5}\;,\\[10pt]
g_{1,12} = \displaystyle \frac{20J_{11}^2 - 16J_{11}J_{22} -
4J_{11}J_{33} - 4J_{22}^2 + 4J_{22}J_{33} -J_{33}^2}{5J_{11}l^3} \;,  \\[10pt]
g_{1,13} = - \displaystyle \frac{3(7K_{33}l^2-480J_{22}
+220J_{11})}{350l^4} + \frac{9(10K_{33}l^2(J_{11} - J_{22})^2 +
112J_{11}J_{22}J_{33})}{35J_{33}K_{33}l^6} \;,  \\[10pt]
g_{1,14} = \displaystyle \frac{12(J_{11}-J_{22})}{l^4} \;,  \\[10pt]
g_{1,15} = - \displaystyle \frac{7K_{33}l^2 + 900( J_{11}+ J_{22})
-700J_{33}}{700l^3} + \frac{ 118K_{33}l^2(J_{11}-J_{22})^2 +
1428J_{11}J_{22}J_{33}}{35J_{33}K_{33}l^5} \;,  \\[10pt]
g_{1,16} = \displaystyle \frac{K_{33}^2l^4 -
8400J_{22}^2}{1400J_{22}l^4}\;,  \\[10pt]
g_{1,17} = - \displaystyle \frac{5J_{11}K_{33}l^2 -
2J_{22}K_{33}l^2 + J_{33}K_{33}l^2
-240J_{11}^2 + 60J_{11}J_{22}+60J_{11}J_{33}}{60J_{11}l^3}\;,  \\[10pt]
g_{1,18} = - \displaystyle \frac{ 7K_{33}l^2 -240J_{11}
-30J_{22}}{1050l^2} + \frac{ 40K_{33}l^2(J_{11}-J_{22})^2 +
462J_{11}J_{22}J_{33}}{35J_{33}K_{33}l^4}\;,  \\[10pt]
g_{1,19} = - \displaystyle \frac{ 7K_{33}l^4 - 270J_{22}K_{33}l^2
- 13860J_{22}^2}{1050K_{33}l^4}\;,  \\[10pt]
g_{1,20} = - \displaystyle \frac{ 10J_{11}^2 -16J_{11}J_{22}
+J_{11}J_{33} -4J_{22}^2+4J_{22}J_{33}-
J_{33}^2}{10J_{11}l^2};\\[10pt]
g_{2,7} = - \displaystyle \frac{ 6(7K_{33}l^2 -480J_{22}
+220J_{11})}{350l^4} + \frac{ 18(10K_{33}l^2(J_{11}-J_{22})^2
+112J_{11}J_{22}J_{33})}{35J_{33}K_{33}l^6}\;,  \\[10pt]
g_{2,10} = \displaystyle \frac{18(7K_{33}^2l^4 -160J_{11}K_{33}l^2
-560J_{11}^2)}{175K_{33}l^7} \;,  \\[10pt]
g_{2,11} = \displaystyle \frac{ 9(7K_{33}^2l^4-260J_{11}K_{33}l^2
-3360J_{11}^2)}{350K_{33}l^6}  \;,  \\[10pt]
g_{2,12} = - \displaystyle \frac{ K_{33}^2l^4 + 840J_{11}K_{33}
-25200J_{11}^2}{700J_{11}l^5} \;,  \\[10pt]
g_{2,13} = \displaystyle \frac{ 63K_{33}^2l^4-520J_{11}K_{33}l^2
-38640J_{11}^2}{700K_{33}l^5} \;,  \\[10pt]
g_{2,14} = \displaystyle \frac{ 14K_{33}l^2 -
 500J_{22}- 80J_{11} + 175 J_{33}}{175J_{33}K_{33}l^5} -
\frac{52K_{33}l^2(J_{11}-J_{22})^2
+ 504J_{11}J_{22}J_{33}}{35J_{33}K_{33}l^5} \;,  \\[10pt]
g_{2,15} = \displaystyle \frac{ 20J_{22}^2
-16J_{11}J_{22}-4J_{22}J_{33} -4J_{11}^2+4J_{11}J_{33} -J_{33}^2
}{5J_{22}l^3} \;,\\[10pt]
g_{2,16} = - \displaystyle \frac{ K_{33}^2l^4 -
8400J_{11}^2}{1400J_{11}l^4} \;,
\end{array} \nonumber
\end{align}

\begin{align} \begin{array}{l}
g_{2,17} = - \displaystyle \frac{ 5J_{22}K_{33}l^2
-2J_{11}K_{33}l^2 +J_{33}K_{33}l^2 -240J_{22}^2 + 60J_{11}J_{22}
+60J_{22}J_{33}}{60J_{22}l^3} \;,  \\[8pt]
g_{2,18} = \displaystyle \frac{7K_{33}l^4
-270J_{11}K_{33}l^2-13860J_{11}^2}{1050K_{33}l^4} \;,  \\[8pt]
g_{2,19} = \displaystyle \frac{ 7K_{33}l^2 - 240J_{22}
-30J_{11}}{1050l^2} - \frac{40K_{33}l^2(J_{11}-J_{22})^2 +
462J_{11}J_{22}J_{33}}{35 l^2} \;, \\[8pt]
g_{2,20} = \displaystyle
\frac{10J_{22}^2-16J_{11}J_{22}+J_{22}J_{33}-4J_{11}^2+4J_{11}J_{33}
-J_{33}^2}{10J_{22}l^2} \;,  \\[8pt]
g_{3,14} = - \displaystyle \frac{ K_{33}(11K_{33}l^2
- 840J_{11})}{6300J_{11}l} \;,\\[10pt]
g_{3,15} = - \displaystyle \frac{ K_{33}(11K_{33}l^2
- 840J_{22})}{6300J_{22}l} \;,\\[8pt]
g_{3,16} = \displaystyle \frac{ K_{33}(2J_{11}^2-J_{11}J_{33}
-2J_{22}^2 +J_{22}J_{33})}{120J_{11}J_{22}} \;,  \\[8pt]
g_{4,18} = \displaystyle \frac{ 7K_{33}l^4-180J_{11}K_{33}l^2-7560J_{11}^2}
{1575K_{33}l^3}  \;,  \\[8pt]
g_{4,19} = \displaystyle \frac{ 14K_{33}l^2 -180( J_{11} + J_{22})
+ 175J_{33}}{1575 l} - \frac{285K_{33}l^2(J_{11}-J_{22})^2
 + 3024J_{11}J_{22}J_{33}}{315J_{33}K_{33}l^3} \;,  \\[8pt]
g_{4,18} = \displaystyle \frac{
12J_{11}^2+28J_{11}J_{22}-12J_{11}J_{33}
-20J_{22}^2+2J_{22}J_{33}+3J_{33}^2}{60J_{22}l}  \;,  \\[8pt]
g_{5,19} = - \displaystyle
\frac{7K_{33}l^4-180J_{22}K_{33}l^2-7560J_{22}^2}
{1575K_{33}l^3} \;,  \\[8pt]
g_{5,20} = - \displaystyle \frac{
12J_{22}^2+28J_{11}J_{22}-12J_{22}J_{33}-20J_{11}^2
+2J_{11}J_{33}+3J_{33}^2}{60J_{11}l}\;,
\end{array} \nonumber \end{align}
and
\begin{align} \begin{array}{lllll}
g_{2,5} = g_{1,7},\; \; &  g_{2,6} = \frac{1}{2} g_{1,8},\; \; &
g_{2,8} = g_{1,14},\; \; & g_{2,9} = g_{1,15},\; \; & g_{3,7} =
g_{1,9}, \\[4pt]
g_{3,8} = g_{1,14},\; \; &  g_{3,9} = \frac{1}{2} g_{1,16},\; \; &
g_{3,10} = g_{1,17},\; \; & g_{3,11} = g_{2,12},\; \; & g_{3,12} =
2 g_{1,16}, \\[4pt]
g_{3,13} = - g_{2,17},\; \; &  g_{4,5} = \frac{1}{2} g_{1,8},\; \;
& g_{4,6} = g_{1,10},\; \; & g_{4,7} = g_{1,15},\; \; & g_{4,8} =
g_{1,17}, \\[4pt]
g_{4,9} = 3 g_{1,18},\; \; &  g_{4,10} = -\frac{1}{3}g_{2,11},\;
\; & g_{4,11} = g_{2,13},\; \; & g_{4,12}
= g_{2,16},\; \; & g_{4,13} = 3 g_{2,18}, \\[4pt]
g_{4,14} = g_{2,19},\; \; &  g_{4,15} = g_{2,20},\; \; & g_{4,16}
= g_{3,14},\; \; & g_{4,17} = g_{3,16},\; \; & g_{5,5} =
\frac{1}{3}g_{1,6}, \\[4pt]
g_{5,6} = - g_{1,11},\; \; &  g_{5,7} = - \frac{1}{2}g_{2,7},\; \;
& g_{5,8} = - g_{1,15},\; \; & g_{5,9}
= - g_{1,16},\; \; & g_{5,10} = - g_{1,18}, \\[4pt]
g_{5,11} = -3 g_{1,19},\; \; &  g_{5,12} = - g_{1,20},\; \; &
g_{5,13} = - g_{2,14},\; \; & g_{5,14}
= - g_{2,17},\; \; & g_{5,15} = - g_{2,19}, \\[4pt]
g_{5,16} = g_{3,15},\; \; &  g_{5,17} = - g_{3,16},\; \; &
g_{5,18} = - g_{4,19},\; \; & g_{6,5}
= g_{1,12},\; \; & g_{6,6} = g_{1,14}, \\[4pt]
g_{6,7} = g_{1,17},\; \; &  g_{6,8} = 2 g_{1,20},\; \; & g_{6,9} =
g_{2,15},\; \; & g_{6,10}
= - g_{3,13},\; \; & g_{6,11} = 2 g_{2,20}, \\[4pt]
g_{6,12} = g_{3,16},\; \; &  g_{6,13} = - g_{4,20},\; \; &
g_{6,14} = g_{5,20}.\; \; &  &
\end{array} \nonumber \end{align}

\renewcommand{\baselinestretch}{1.2}

\end{document}